\documentclass[oneside,final]{csri23}

 \usepackage[top=1.355in,bottom=1.355in,left=1.5in,right=1.5in]{geometry}
\usepackage{amsfonts,amsmath,graphicx,url}
\usepackage{float}
\usepackage{color}
\usepackage[T1]{fontenc}
\usepackage{ae,aecompl}
\usepackage[abs]{overpic}
\usepackage[normalem]{ulem}
\usepackage[font=footnotesize]{caption}
\usepackage{subcaption}

\usepackage[ruled,vlined]{algorithm2e}

\title{Domain Decomposition-based Coupling of Physics-Informed Neural Networks via the Schwarz Alternating Method}

\author{William Snyder\thanks{Virginia Polytechnic Institute, swilli9@vt.edu} \and Irina Tezaur\thanks{Sandia National Laboratories,
ikalash@sandia.gov} \and Christopher R. Wentland\thanks{Sandia National Laboratories, crwentl@sandia.gov}}

\begin{document}

\maketitle

\begin{abstract}
Physics-informed neural networks (PINNs) are appealing data-driven tools for solving and inferring solutions to nonlinear partial differential equations (PDEs). Unlike traditional neural networks (NNs), which train only on solution data, a PINN incorporates a PDE's residual into its loss function and trains to minimize the said residual at a set of collocation points in the solution domain. This paper explores the use of the Schwarz alternating method as a means to couple PINNs with each other and with conventional numerical models (i.e., full order models, or FOMs, obtained via the finite element, finite difference or finite volume methods) following a decomposition of the physical domain. It is well-known that training a PINN can be difficult when the PDE solution has steep gradients.
We investigate herein the use of domain decomposition and the Schwarz alternating method as a means to accelerate the PINN training phase. Within this context, we explore different approaches for imposing Dirichlet boundary conditions within each subdomain PINN: weakly through the loss and/or strongly through a solution transformation.  As a numerical example, we consider the one-dimensional steady state advection-diffusion equation in the advection-dominated (high P\'{e}clet) regime.  Our results suggest that the convergence of the Schwarz method is strongly linked to the choice of boundary condition implementation within the PINNs being coupled.  Surprisingly, strong enforcement of the Schwarz boundary conditions does not always lead to a faster convergence of the method.  While it is not clear from our preliminary study that the PINN-PINN coupling via the Schwarz alternating method accelerates PINN convergence in the advection-dominated regime, it reveals that PINN training \textit{can} be improved substantially for P\'{e}clet numbers as high as $1.0\times 10^6$ by performing a PINN-FOM coupling.
\end{abstract}

\section{Introduction} \label{WDS:sec:intro}

In recent years, physics-informed neural networks (PINNs)~\cite{WDS:Raissi2019} have emerged as a scientific machine learning (ML) technique advertised to solve partial differential equations (PDEs) without requiring any training data or an underlying mesh. This is accomplished by constructing and minimizing a loss function that includes the residual of the underlying PDE, sampled at a finite number of collocation points within the physical domain on which the PDE is posed. While the aforementioned properties of PINNs make them appealing within the modeling and simulation community, these models unfortunately suffer from several deficiencies.  First, PINNs are notoriously difficult to train for problems having a slowly-decaying Kolmogorov $n$-width, as discussed in~\cite{WDS:Mojgani:2023} and the references therein.  This class of problems includes advection-dominated flow problems in which the solution exhibits sharp boundary layers and/or shocks. Second, since PINNs are trained for a given geometry with a given set of boundary conditions, they require retraining when applied to different geometries or boundary data~\cite{WDS:Wang:2022}.

In this paper, we present and evaluate a promising technique that can mitigate both difficulties described above: the Schwarz alternating method~\cite{WDS:schwarz1870ueber} for domain decomposition-based coupling. The Schwarz alternating method is based on the simple idea that if the solution to a PDE is known in two or more regularly-shaped subdomains comprising a more complex domain, these local solutions can be used to iteratively build a solution on the more complex domain, with information propagating between subdomains through  boundary conditions imposed on the subdomain boundaries. In several of our recent works~\cite{WDS:mota2017schwarz, WDS:mota2022schwarz}, the overlapping version of the Schwarz alternating method was pioneered as a mechanism for performing concurrent coupling of subdomains discretized in space by disparate meshes and in time by different time-integration schemes with different time-steps. In those works, attention was restricted to the coupling of high-fidelity full order models (FOMs) in quasistatic and dynamic solid mechanics. Our recent work~\cite{WDS:barnett2022schwarz} extended the approach to the coupling of projection-based reduced order models (ROMs) with each other and with FOMs following either an overlapping or a non-overlapping domain decomposition (DD) of the underlying geometry. We additionally demonstrated in~\cite{WDS:Koliesnikova2023, WDS:Hoy2021} that it is possible to transform the non-overlapping Schwarz alternating method into a novel contact enforcement algorithm which exhibits remarkable energy conservation properties.  

The present work explores the use of the Schwarz alternating method for coupling PINNs with each other and with FOMs, following an overlapping DD of the underlying physical domain. For a detailed literature overview on existing methods which combine ML and domain decomposition methods, the reader is referred to~\cite{WDS:Heinlein:2020}. While it is possible to couple pre-trained subdomain-local PINNs or NNs, as in~\cite{WDS:Wang:2022}, our initial study focuses primarily on evaluating the Schwarz alternating method as a means to facilitate PINN training.  The proposed approach is similar to the deep domain decomposition method (D3M)  and the deep domain decomposition (DeepDDM) methods, proposed recently in~\cite{WDS:LiD3M} and~\cite{WDS:LiDeepDDM}, respectively. Unlike the D3M method, which is based on a deep Ritz method PINN formulation in which the loss function is derived from the weak variational form of the governing PDEs~\cite{WDS:e2017deep, WDS:Sukumar:2022}, our loss function incorporates the PDE residual in strong form. Unlike the DeepDDM method, which imposes Schwarz and other boundary conditions weakly through a boundary loss contribution to the loss function being minimized, we consider both strong and weak formulations of the Dirichlet boundary conditions (DBCs) within our algorithm, as well as a combination of the two. For the former strong DBC enforcement, we borrow ideas from the Finite Basis PINN (FBPINN) literature~\cite{WDS:Moseley:2023, WDS:Dolean:2023}. While the authors of~\cite{WDS:LiD3M, WDS:LiDeepDDM} focus their development and demonstrations on the Poisson equation, our work considers a much more difficult boundary value problem (BVP), namely the advection-diffusion equation in the advection-dominated (high P\'{e}clet) regime.  Finally, to the best of our knowledge, ours is the first coupling formulation that enabling the creation of ``hybrid'' models through the coupling of PINNs with conventional numerical models, such as those constructed using the finite element, finite difference, or finite volume method. 

The remainder of this paper is organized as follows. Section~\ref{WDS:sec:addiff} defines the model problem considered herein, a one-dimensional (1D) advection-diffusion equation on the unit interval. Section~\ref{WDS:sec:Schwarz} summarizes the overlapping version of the Schwarz alternating method for domain decomposition-based coupling in the specific context of our advection-diffusion model problem. Section~\ref{WDS:sec:pinns} provides an overview of PINNs and several PINN variants, namely PINNs in which DBCs are enforced strongly, rather than weakly, and PINNs in which a data loss term is incorporated into the loss function being minimized.  In Section~\ref{WDS:sec:pinns_schwarz}, we extend the general Schwarz formulation described in Section~\ref{WDS:sec:Schwarz} to the specific case of PINN-PINN and PINN-FOM coupling.  The proposed method is evaluated as a mechanism for facilitating PINN training via domain decomposition-based coupling on our model advection-diffusion problem in Section~\ref{WDS:sec_numerical}.  
Whereas a definitive conclusion on whether or not
PINN training can be facilitated through 
PINN-PINN coupling using the Schwarz alternating method cannot be made from the present study, our numerical experiments reveal that 
this goal \textit{can} be achieved by performing a Schwarz-based PINN-FOM coupling.
Finally, we provide some conclusions and directions for future work in Section~\ref{WDS:sec:conc}.

\section{Model problem}\label{WDS:sec:addiff}

\begin{figure}[htbp!]
    \begin{minipage}{0.45\linewidth}
        \includegraphics[width=0.99\textwidth]{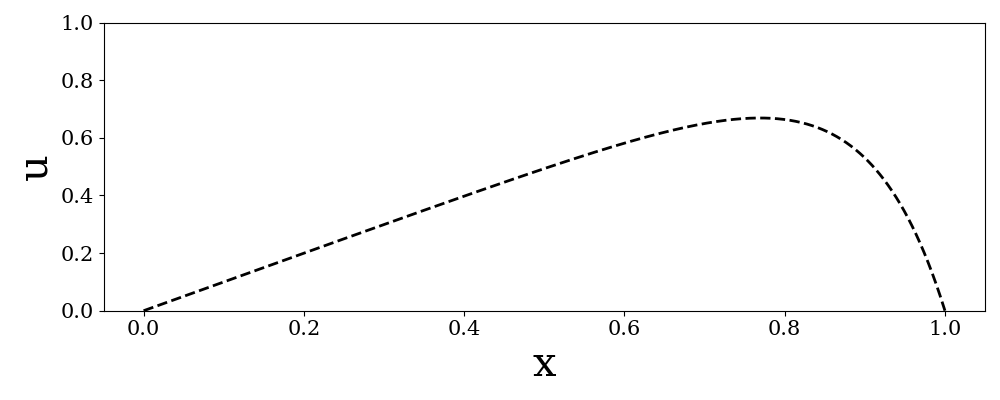}
        \subcaption{$Pe = 10$}
    \end{minipage}
    \begin{minipage}{0.45\linewidth}
        \includegraphics[width=0.99\textwidth]{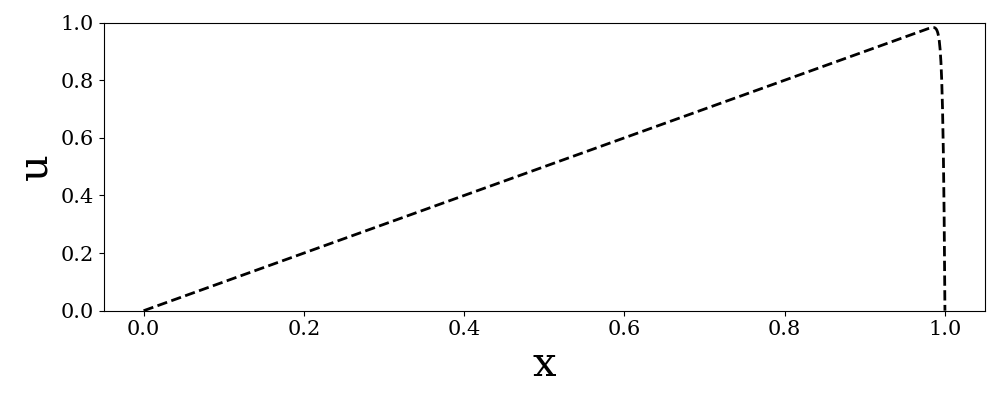}
        \subcaption{$Pe = 500$}
    \end{minipage}
    \begin{center}
    \begin{minipage}{0.45\linewidth}
        \includegraphics[width=0.99\textwidth]{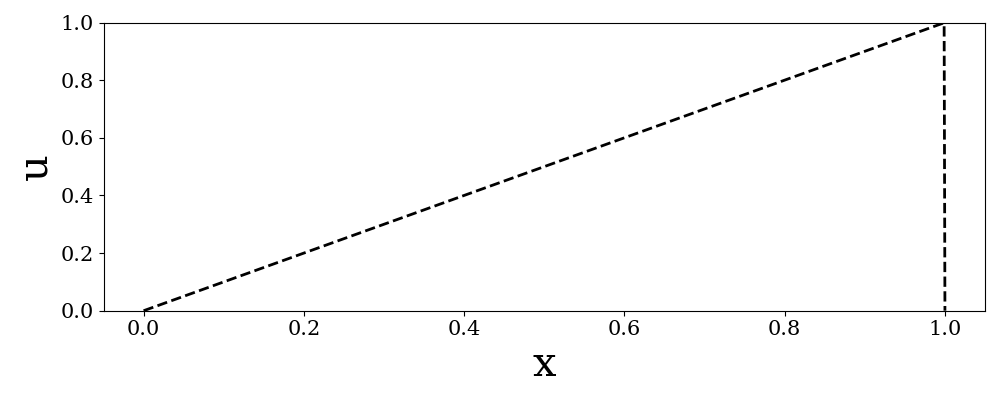}
        \subcaption{$Pe = 1 \times 10^6$}
    \end{minipage}
    \end{center}
	\caption{Solutions to our model advection-diffusion problem for different P\'{e}clet numbers $Pe$.}
    \label{fig:adv-diff-Pe}
\end{figure}

As a numerical example to test our coupling method and PINNs, we choose the 1D steady-state advection-diffusion equation defined in an open bounded domain $\Omega = (0,1)$ 
with boundary $\partial \Omega = \{ 0, 1\}$:
\begin{equation}
	-{\nu}\frac{\partial^2 u}{\partial x^2}+\frac{\partial u}{\partial x}-1=0 \hspace{0.5cm} \text{in } \Omega = (0,1),
    \label{eqn:residual}
\end{equation}
with boundary conditions
\begin{equation}
	u(0)=u(1)=0.
    \label{eqn:BC}
\end{equation}
From this point forward, we will refer to the DBCs in~\eqref{eqn:BC} as the ``system DBCs'' or ``system boundary conditions''. In~\eqref{eqn:residual}, $\nu$ is the diffusion coefficient, which defines the P\'{e}clet number of the problem, given by $Pe:=\nu^{-1}$. Our focus here is primarily on advection-dominated regimes for this problem, i.e., $10 \leq Pe \leq 10^{6}$. Figure~\ref{fig:adv-diff-Pe} shows a plot of the solution to~\eqref{eqn:residual}--\eqref{eqn:BC} as a function of the P\'{e}clet number. The reader can observe that the solution develops a sharp gradient near the right boundary of the domain, whose steepness is a function of $Pe$.

\section{The Schwarz alternating method for domain decomposition-based coupling}\label{WDS:sec:Schwarz}

The Schwarz alternating method is the oldest known domain decomposition method, first proposed in 1870 by H. Schwarz~\cite{WDS:schwarz1870ueber}. The method is based on the simple idea that the solution to a PDE on a complex domain $\Omega$ can be obtained via an iterative procedure on subdomains comprising $\Omega$, with information propagating via transmission (or boundary) conditions imposed on subdomain boundaries. While the Schwarz alternating method can be formulated for both overlapping~\cite{WDS:mota2017schwarz, WDS:mota2022schwarz, WDS:Lions1988} and non-overlapping~\cite{WDS:barnett2022schwarz, WDS:Koliesnikova2023, WDS:Hoy2021} DDs, we restrict our attention to the overlapping variant of the method, which has been shown to have favorable convergence properties with respect to its non-overlapping counterpart \cite{WDS:barnett2022schwarz}.

\begin{figure}[htbp!]
        \begin{minipage}{0.49\linewidth}
            \includegraphics[width=0.99\textwidth]{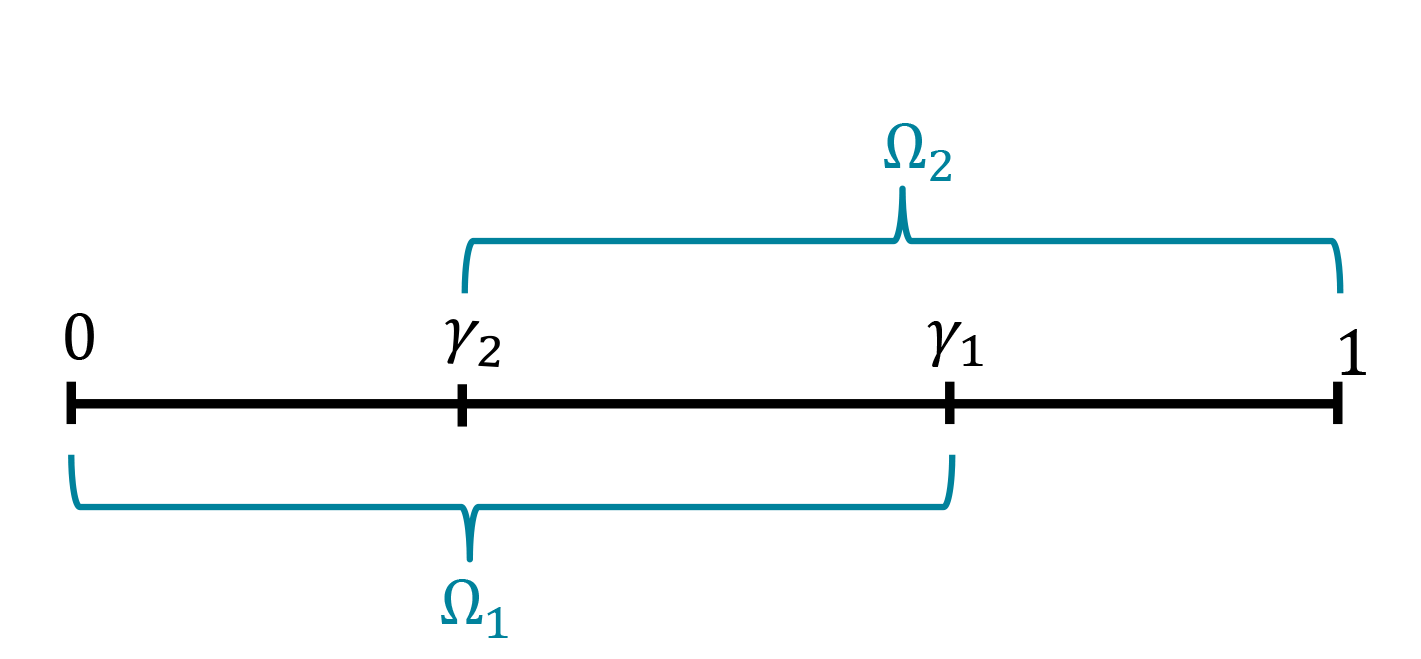}
            \subcaption{1D}
        \end{minipage}
        \begin{minipage}{0.49\linewidth}
            \includegraphics[width=0.99\textwidth]{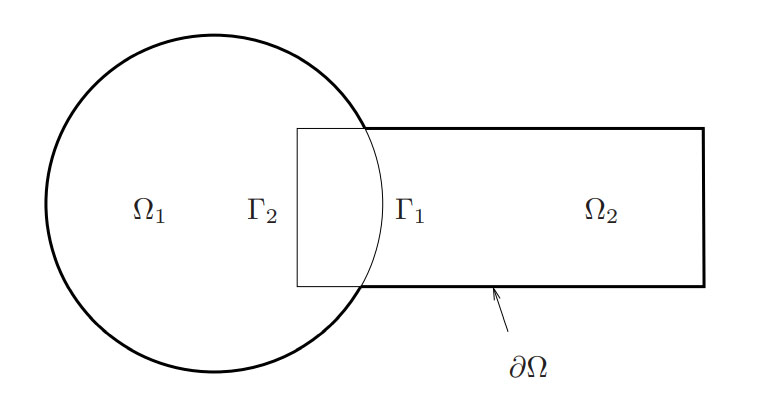}
            \subcaption{2D}
        \end{minipage}
            
	\caption{Example overlapping domain decompositions for Schwarz alternating method in 1D (left) and 2D (right).}
    \label{fig:dds}
\end{figure}

Consider 
an overlapping DD of a geometry $\Omega$ into two subdomains, as shown in Figures~\ref{fig:dds}(a) and (b) in one and two spatial dimensions, respectively. We present the basic overlapping Schwarz algorithm on our targeted model problem~\eqref{eqn:residual} with system boundary conditions~\eqref{eqn:BC}.  Suppose the governing domain for our model problem is partitioned into two overlapping domains: $\Omega = \Omega_1 \cup \Omega_2$, with $\Omega_1 = (0,\gamma_1)$ and $\Omega_2 = (\gamma_2, 1)$, with $0 < \gamma_2 < \gamma_1 < 1$, so that $\Omega_1 \cap \Omega_2 \neq \emptyset$, as shown in Figure \ref{fig:dds}(a). Equations~\eqref{eqn:residual}--\eqref{eqn:BC} can be solved using the Schwarz alternating method by performing the following iteration:

\begin{equation}\label{eq:advdiff_schwarz}
    \left\{
    \begin{aligned}
	    -\nu \frac{\partial^2 u_1^{(n+1)}}{\partial x^2} + \frac{\partial u_1^{(n+1)}}{\partial x}  & = 1, \text{ in } \Omega_1,
        \\
        u_1^{(n+1)} & = 0,  \text{ at } x = 0,
        \\
        u_1^{(n+1)} & = u_2^{(n)},  \text{ at }  x=\gamma_1,
        \\
    \end{aligned}
    \right. \hspace{0.2cm}
    \left\{
    \begin{aligned}
	    -\nu \frac{\partial^2 u_2^{(n+1)}}{\partial x^2} + \frac{\partial u_2^{(n+1)}}{\partial x}  & = 1, \text{ in } \Omega_2,
        \\
        u_2^{(n+1)} & = 0,  \text{ at } x = 1,
        \\
        u_2^{(n+1)} & = u_1^{(n+1)},  \text{at }  x = \gamma_2, 
        \\
    \end{aligned} \right.
\end{equation}
for Schwarz iterations $n = 0, 1, 2, ...$.  In~\eqref{eq:advdiff_schwarz}, $u_i$ denotes the solution in $\Omega_i$ for $i=1,2$.  The reader can observe that the transmission boundary conditions on the Schwarz boundaries $x = \gamma_1$ and $x = \gamma_2$ are of the Dirichlet type. The iteration~\eqref{eq:advdiff_schwarz} continues until convergence is reached. Generally, convergence of the method is declared\footnote{Herein, for the data-driven models coupled via the Schwarz alternating method, 
we also require the solution to be within a given tolerance from the known analytical or high-fidelity solution to the given problem; more details are provided in Section \ref{WDS:sec_numerical}.} when $||u_i^{(n+1)} - u_i^{(n)}||_2/||u_i^{(n)}||_2 < \delta$, for $i=1,2$ and for some specified tolerance $\delta$. In the remainder of this paper, the boundaries $\gamma_i$ will be referred to as the ``Schwarz boundaries" and the boundary conditions applied at these boundaries will be referred to as the ``Schwarz boundary conditions".  
It can be shown~\cite{WDS:Lions1988, WDS:mota2017schwarz} that the Schwarz iteration procedure~\eqref{eq:advdiff_schwarz} converges to the solution of~\eqref{eqn:residual}--\eqref{eqn:BC} provided the overlap is non-empty ($\Omega_1 \cap \Omega_2 \neq \emptyset$). It is straightforward to extend the Schwarz iteration procedure~\eqref{eq:advdiff_schwarz} to an arbitrary number of overlapping subdomains. The present work builds on recent extensions of the Schwarz alternating method to concurrent multi-scale coupling in quasistatic~\cite{WDS:mota2017schwarz} and dynamic~\cite{WDS:mota2022schwarz} solid mechanics, and to the coupling of projection-based reduced order models (ROMs) with each other and with high-fidelity full order models (FOMs)~\cite{WDS:barnett2022schwarz}.  \\

\noindent {\bf \textit{Remark 1}}.  
The Schwarz iteration procedure described above in~\eqref{eq:advdiff_schwarz} is often referred to as the multiplicative Schwarz algorithm~\cite{WDS:Gander2008}. In this algorithm, the solution in subdomain $\Omega_2$ at iteration $n+1$ depends on the solution in subdomain $\Omega_1$ also at iteration $n+1$. This contemporaneous dependence prevents one from solving the $\Omega_1$ and $\Omega_2$ problems in parallel.  The Schwarz iteration can be modified to achieve what is commonly referred to as additive Schwarz~\cite{WDS:Gander2008} by applying boundary conditions from the $n^{th}$ Schwarz iteration procedure in $\Omega_1$ to the $(n + 1)^{st}$ $\Omega_2$ sub-problem. If this change is made, the Schwarz iteration sequence within can be parallelized over the number of subdomains, as done in~\cite{WDS:LiD3M, WDS:LiDeepDDM}. While we do not consider the additive variant of the Schwarz method in the present work, preliminary results have suggested that the method does not reduce solution accuracy and can achieve speed-ups if parallelized appropriately.

\section{Physics-informed neural networks (PINNs)} \label{WDS:sec:pinns}

\begin{figure}[htbp!]
    \begin{center}       
        \includegraphics[width=0.96\textwidth]{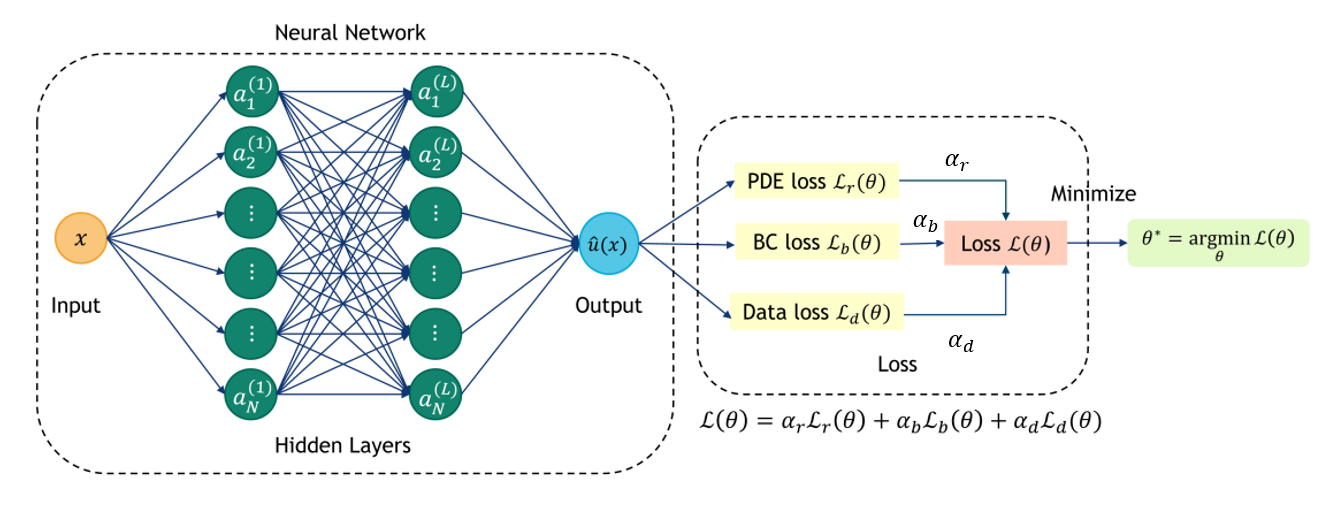}
    \end{center}
	\caption{Illustration of our PINN architecture for $L = 2$ layers each having $N = 6$ nodes per layer.}
    \label{fig:PINN}
\end{figure}

PINNs, introduced in the seminal paper~\cite{WDS:Raissi2019}, are constructed via a least-squares collocation method applied to the governing PDEs, the solution of which is approximated by a neural network function, denoted by $NN(x; \theta)$. Let $x$ denote the input and $\hat{u}$ denote the output of a typical feed-forward fully-connected NN\footnote{From this point further, we will use the term NN and PINN interchangeably, since our most general PINN formulation~\eqref{eq:loss} includes a data loss term.  Please see Remark 2 for a comment on this.} with $L$ hidden layers and $N$ neurons in each layer, as shown in Figure~\ref{fig:PINN}. The output of the network is computed as
\begin{equation} \label{eq:weights_biases1}
    NN(x ; \theta) = \sigma \left(\sum_{i=1}^N w_{1,i}^{(L)} a_i^{(L)} + b^{(L)} \right), 
\end{equation}
where
\begin{equation} \label{eq:weights_biases2}
    a_n^{(j)} = \sigma\left(\sum_{i=1}^N w_{k,i}^{(j-1)} a_i^{(j-1)} + b^{(j-1)} \right),
\end{equation}
for $j = 2, ..., L$ and $k = 1, .., N$, and 
\begin{equation} \label{eq:weights_biases3}
    a_n^{(1)} = \sigma\left(w_{k,1}^{(0)}x + b^{(0)} \right),
\end{equation}
for $k = 1, .., N$. In~\eqref{eq:weights_biases1}--\eqref{eq:weights_biases3}, $\sigma(\cdot)$ denotes a nonlinear function known as the activation function, $a_i^{(j)}$ is referred to as the activation of neuron $i$ in layer $j$, and the parameters $\theta := \left( w_{1,1}^{(0)}, ..., w_{1,N}^{(L)}, b^{(0)}, ..., b^{(L)} \right)^T$ include the weights ($w_{i,j}^{(l)}$) and biases ($b^{(l)}$) of the NN. The parameters $\theta$ are computed by minimizing a pre-defined loss function $\mathcal{L}(\theta)$:
\begin{equation} \label{eq:min}
    \theta^* = \arg\min_{\theta} \mathcal{L}(\theta).
\end{equation}
A PINN is differentiated from a typical NN by the definition of the loss function $\mathcal{L}(\theta)$. Whereas a NN loss function measures the difference between the NN outputs and a corresponding set of solution data taken as the ground truth, a PINN does not require any such solution data. Instead, a PINN uses the residual of the PDEs it is attempting to solve as the loss function (in this case, equation~\eqref{eqn:residual}). This is achieved by first taking as inputs a set of collocation points $x\in \Omega$. In its most basic form, the loss function defining the PINN minimization problem is given by the residual loss, that is, by 
\begin{equation} \label{eq:loss_simple}
     \mathcal{L}(\theta) = \mathcal{L}_r(\theta) := \frac{1}{M} \sum_{i=1}^M \left(-{\nu} \nabla_x^2 \hat{u}(x_i,\theta) + \nabla_x \hat{u}(x_i,\theta) - 1 \right)^2, 
\end{equation}
where the $x_i\in \Omega$ for $i=1,...,M$ are the collocation points used to evaluate the residual loss and $NN(x; \theta) = \hat{u}(x; \theta)\approx u(x)$ represents the PINN approximation of the solution as a function of $x$ and the training parameters $\theta$. Once the loss function~\eqref{eq:loss_simple} is defined along with the PINN architecture (Figure~\ref{fig:PINN}), a gradient-based optimization algorithm is used to solve~\eqref{eq:min}, with derivatives computed via automatic differentiation.

In the following subsections, we describe several PINN variations that allow for the incorporation of boundary conditions as well as available snapshot data.  In this more general formulation, the loss being minimized takes the form
\begin{equation} \label{eq:loss}
    \mathcal{L}(\theta) = \alpha_r \mathcal{L}_r(\theta) + \alpha_b \mathcal{L}_b(\theta) 
    + \alpha_d \mathcal{L}_d(\theta),
\end{equation}
for some relaxation parameters $\alpha_r, \alpha_b, \alpha_d \in \mathbb{R}$, often normalized to sum to be between 0 and 1. In~\eqref{eq:loss}, $\mathcal{L}_b(\theta)$ and $\mathcal{L}_d(\theta)$ are the boundary and data losses, respectively. These terms are discussed in more detail below, in Sections \ref{WDS:sec_WDBC} and \ref{WDS:sec_data}, respectively.  \\

\noindent {\bf \textit{Remark 2.}}  While, in the original PINN paper by Raissi \textit{et al.}~\cite{WDS:Raissi2019}, PINNs were presented as completely data-free, i.e., $\alpha_d = 0$ in~\eqref{eq:loss}, subsequent PINN formulations allowed the addition of the data-loss term $\mathcal{L}_d(\theta)$ for cases where some training data (from high-fidelity simulations or experiments) are available. For more information, the interested reader is referred to~\cite{WDS:Lawal2022, WDS:Hao2023physicsinformed} and the references therein. \\

\subsection{Weak Dirichlet BCs (WDBCs)}\label{WDS:sec_WDBC}

The traditional way to impose boundary conditions in PINNs is through a weak formulation, via the loss function. To this effect, the loss function being minimized~\eqref{eq:loss} takes the form 
\begin{equation} \label{eq:loss2}
    \mathcal{L}(\theta) = \alpha \mathcal{L}_r(\theta) + (1-\alpha)\mathcal{L}_b(\theta),
\end{equation}
where $\alpha \in (0,1]$, and 
\begin{equation} \label{eq:boundary_loss}
    \mathcal{L}_b(\theta) := \frac{1}{b}\sum_{i=1}^b\left(\hat{u}(x_i, \theta) - u(x_i) \right)^2.
\end{equation}
In~\eqref{eq:boundary_loss}, the points $x_i \in \partial \Omega$ for $i = 1, ..., b$ are the coordinates of boundary points at which Dirichlet boundary conditions are prescribed. In the 1D domain depicted in Figure~\ref{fig:dds}(a), $x_i \in \{ 0,1\}$ for the domain $\Omega$, whereas $x_i \in \{0, \gamma_1\}$ and $x_i \in \{ \gamma_2, 1\}$ for the subdomains $\Omega_1$ and $\Omega_2$, respectively. For 1D subdomains, the maximum number of boundary points for which the weak boundary loss is computed (within each subdomain) is thus $b = 2$.

\subsection{Strong Dirichlet BCs (SDBCs)} \label{WDS:sec_SDBC}

While the majority of PINNs in the literature employ a weak enforcement of the boundary conditions via a boundary loss term~\eqref{eq:boundary_loss}, this approach has some disadvantages. First, it is unclear \textit{a priori} how to select the relaxation parameters $\alpha_r$ and $\alpha_b$ in~\eqref{eq:loss}; often this is done by trial and error, and the values require retuning when the problem setup is modified.  Second, since the BCs in \eqref{eq:loss} are imposed weakly, the learned solution may be inconsistent with the underlying BVP~\eqref{eqn:residual}--\eqref{eqn:BC}. Moreover, some recent work has demonstrated that the PINN optimization problem~\eqref{eq:loss} may be extremely stiff, leading to a lack of convergence as a result of the residual and boundary losses competing with each other in the loss function~\cite{WDS:Wang:2021, WDS:Sun:2020}.

Since information within the Schwarz alternating method propagates through the DBCs imposed at the subdomain boundaries and past convergence analyses~\cite{WDS:Lions1988, WDS:mota2017schwarz, WDS:mota2022schwarz} of the method have assumed a strong implementation of these BCs, we are interested in being able to implement SDBCs within our PINNs for the BVP considered. 

Motivated by earlier related work for other flavors of 
collocation methods (e.g., \cite{WDS:Graepel:2003} for Bayesian Gaussian processes), recent years have seen the development of several approaches for imposing BCs strongly within a PINN or NN.  In the FBPINN approach~\cite{WDS:Moseley:2023,WDS:Dolean:2023}, instead of defining a NN to directly approximate the solution as discussed above ($NN(x;\theta) \approx u(x)$), an ansatz of the form
\begin{equation} \label{eq:ansatz_fbpinn}
    \hat{u}(x; \theta) = g(x) + \psi(x) NN(x; \theta) \approx u(x),
\end{equation}
where the functions $g(x)$ and $\psi(x)$ are derived such that $\hat{u}(x,\theta)$ in~\eqref{eq:ansatz_fbpinn} satisfies the prescribed DBCs exactly. Since DBCs are imposed strongly in the NN, it is not necessary to include the boundary loss $\mathcal{L}_b(\theta)$ in the loss function; hence, the loss function reduces to the residual loss, i.e.,~\eqref{eq:loss_simple}. In the case of complex, multi-dimensional domains, the functional form of $g(x)$ may be rather complex, and will often depend on the distance of a point $x$ from the boundary $\partial \Omega$, as discussed in~\cite{WDS:Wang:2022}. An alternate approach that uses distance functions and geometry-aware trial functions within a Deep Ritz PINN (i.e., a PINN in which the residual loss is given in weak variational rather than strong form) is the work of Sukumar \textit{et al.} \cite{WDS:Sukumar:2022}.

In order to keep the discussion as general as possible for our 1D model problem~\eqref{eqn:residual}, consider a generic overlapping decomposition of $\Omega$ into $n_{D}$ overlapping subdomains, so that $\Omega = \cup_{i=1}^{n_{D}} \Omega_i$, where $\Omega_i = (\gamma_{2i-2}, \gamma_{2i-1})$, for $i=1, ..., n_D$, where $\gamma_0 = 0 $ and $\gamma_{2n_D-1} = 1$.  We begin by defining the following subdomain-local BVP for our targeted 1D advection-diffusion equation
\begin{equation}
    -{\nu}\frac{\partial^2 u_i}{\partial x^2}+\frac{\partial u_i}{\partial x}-1=0, \hspace{0.5cm} \text{in  }\Omega_i:= (\gamma_{2i-2}, \gamma_{2i-1}),
\end{equation}
with boundary conditions
\begin{equation} \label{eq:DBC_Omegai}
    u_i(\gamma_{2i-2}) = g_{2i-2}, \hspace{0.5cm} u_i(\gamma_{2i-1}) = g_{2i-1},  
\end{equation}
for $g_i \in \mathbb{R}$, with $i=1,..., n_D$.  Let $NN_i(x; \theta)$ denote a NN trained to represent the solution in $\Omega_i$.  The reader can observe that the following function is guaranteed to satisfy strongly the DBCs~\eqref{eq:DBC_Omegai} on $\partial \Omega_i$:
\begin{equation}
    \hat{u}_i (x; \theta) = v_i(x) NN_i(x; \theta) + \phi_i(x) g_{2i-2} + \psi_i(x) g_{2i-1},
\end{equation}
where $v_i(x)$ is a function such that $v(\gamma_{2i-2}) = v_i(\gamma_{2i-1}) = 0$, $\phi_i(x)$ is a function such that $\phi_i(\gamma_{2i-2}) = 1 $ and $\phi_i(\gamma_{2i-1}) = 0$, and $\psi_i(x)$ is a function such that $\psi_i(\gamma_{2i-2}) = 0$ and $\psi_i(\gamma_{2i-1}) = 1$.  In the present work, we employ the following expressions for these functions:
\begin{equation}
    v_i(x) = \tanh(k(\gamma_{2i-1}-x))\tanh(k(x - \gamma_{2i-2})),
    \label{eqn:v_x}
\end{equation}
\begin{equation} \label{eq:phi_i}
    \phi_i(x) = 10^{-10(x-\gamma_{2i-2})},
\end{equation}
\begin{equation} \label{eq:psi_i}
    \psi_i(x) = 10^{10(x-\gamma_{2i-1})},
\end{equation}
for $x \in (\gamma_{2i-2}, \gamma_{2i-1})$ and $k>0$. 

\begin{figure}
    \begin{minipage}{0.49\linewidth}
        \includegraphics[width=0.99\textwidth]{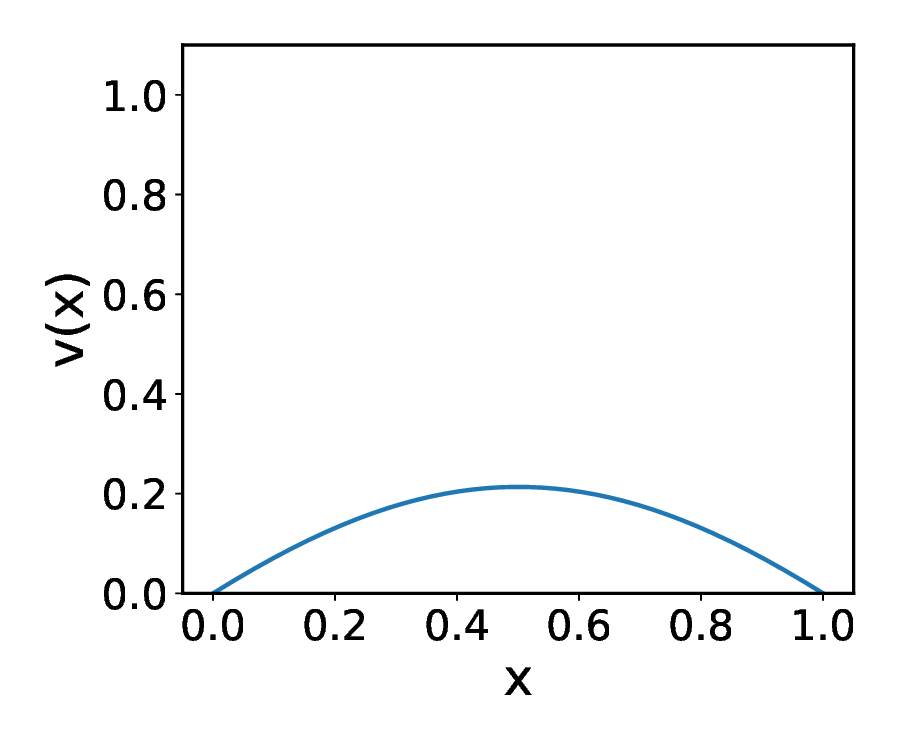}
        \subcaption{$k=1$}
    \end{minipage}
    \begin{minipage}{0.49\linewidth}
        \includegraphics[width=0.99\textwidth]{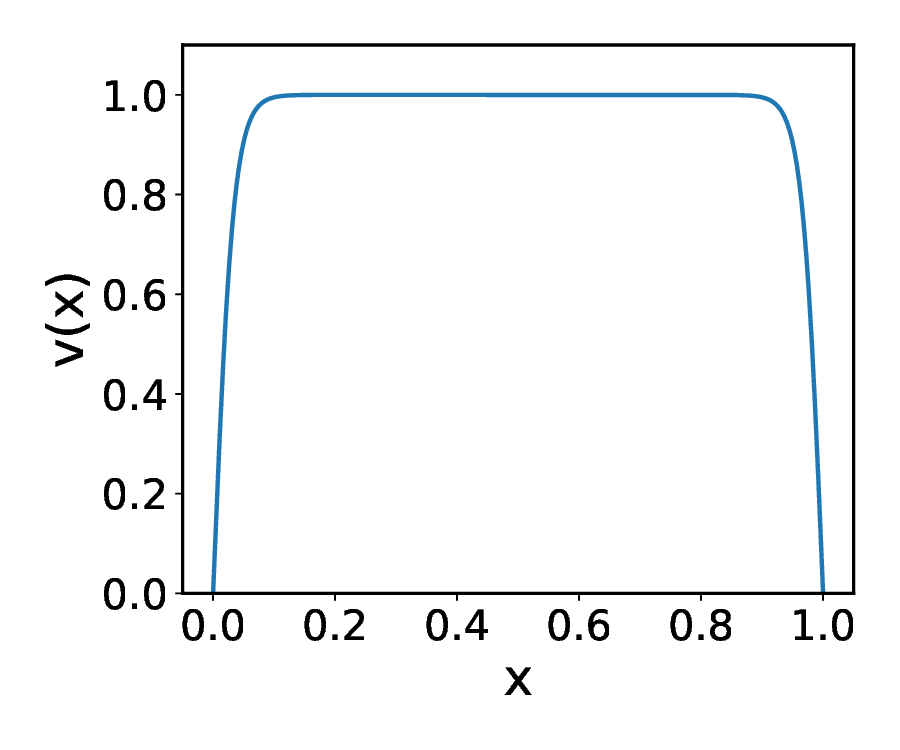}
        \subcaption{$k=30$}
    \end{minipage}
    \caption{Strong DBC scaling function \eqref{eqn:v_x} for two values of $k$: $k=1$ (left) and $k=30$ (right).}
    \label{fig:tanh}
\end{figure}

The scaling function~\eqref{eqn:v_x} is plotted for two values of $k$, $k=1$ and $k=30$, in Figure~\ref{fig:tanh}. As $k$ is increased, 
$v_i(x)$ begins to resemble a step function on $\Omega = (0,1)$.  
At first glance, it may seem as though employing a larger $k$ in~\eqref{eqn:v_x} should improve PINN convergence, as $v_i(x) \approx 1$ in $\Omega$ away from the boundary $\partial \Omega$, meaning that this scaling function minimally modifies the NN away from $\partial \Omega$ (Figure \ref{fig:tanh}(a)). However, our numerical experiments reveal that selecting higher values of $k$, e.g., $k=30$ (Figure~\ref{fig:tanh}(b)), severely hinders the NN's convergence. We believe that this happens because $v_i(x)$ exhibits sharp gradients when $k \gg 1$, requiring the NN to make rapid and substantial changes to the magnitudes of its outputs for a small subset of the spatial domain close to the boundaries, complicating the minimization of the loss function. For this reason, for the results presented in Section~\ref{WDS:sec_numerical}, a value of $k=1$ is used to define $v_i(x)$~\eqref{eqn:v_x}.   

\begin{figure}
    \begin{minipage}{0.49\linewidth}
        \includegraphics[width=0.99\textwidth]{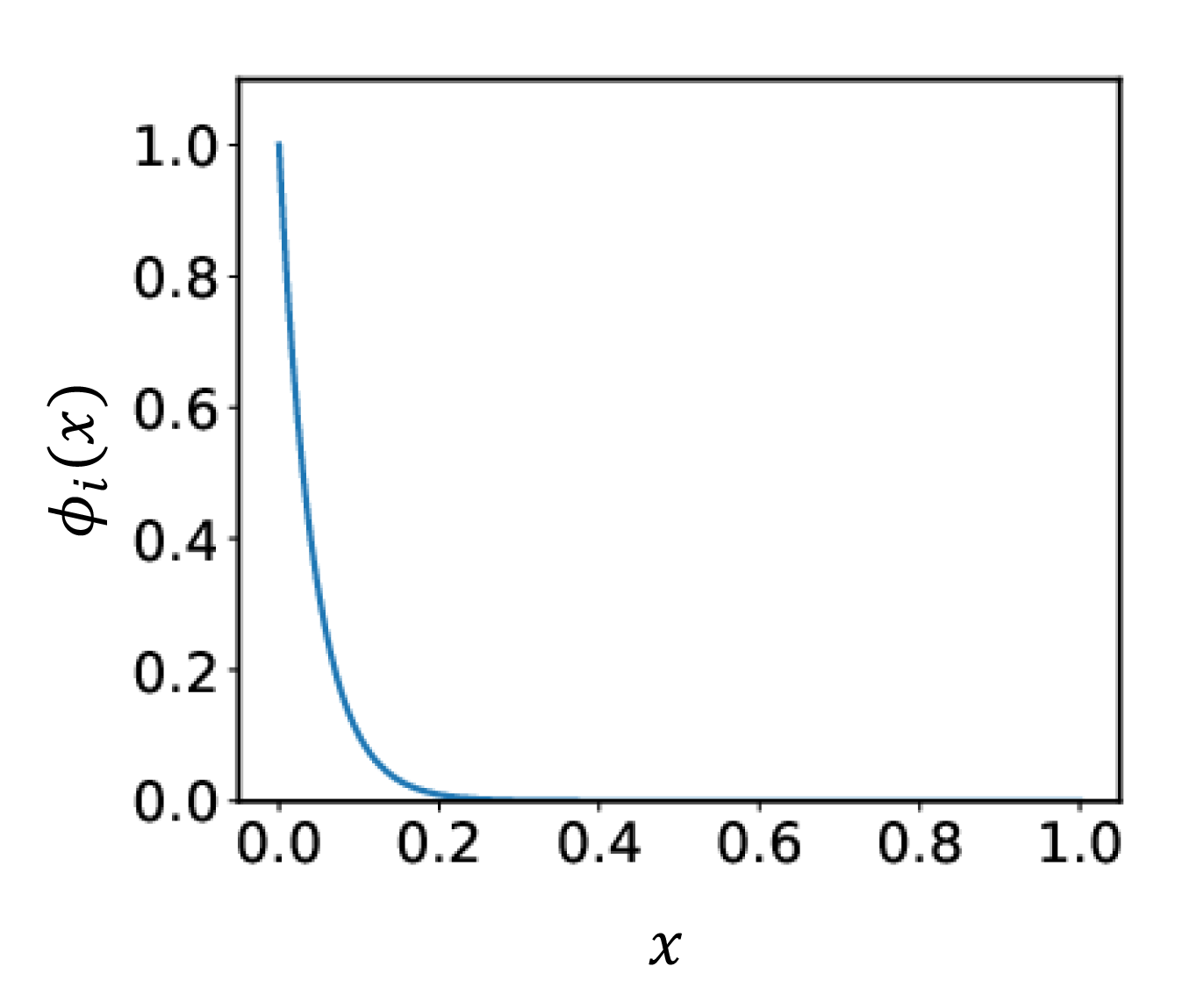}
        \subcaption{$\phi_i(x)$}
    \end{minipage}
    \begin{minipage}{0.49\linewidth}
        \includegraphics[width=0.99\textwidth]{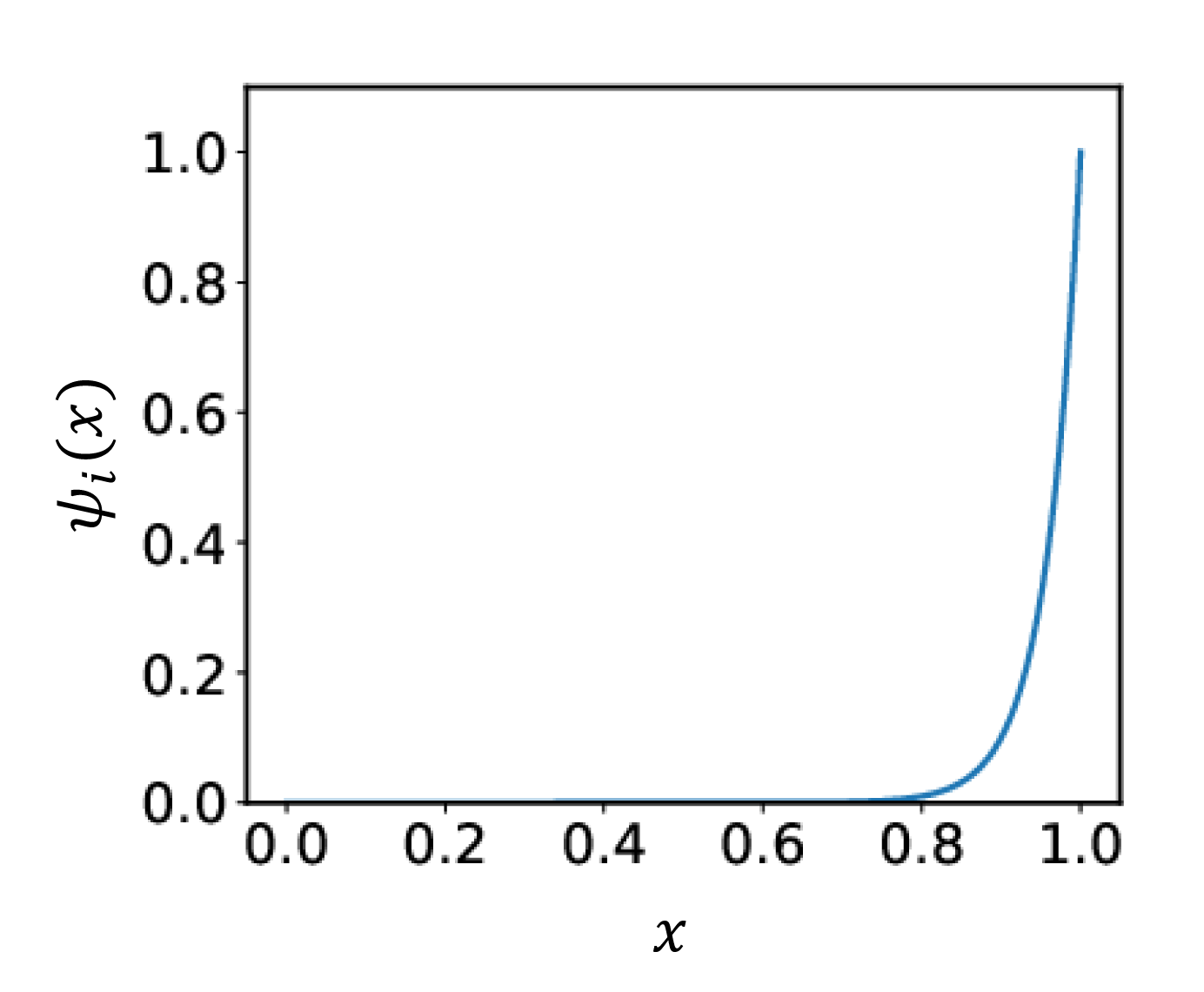}
        \subcaption{$\psi_i(x)$}
    \end{minipage}
    \caption{Inhomogeneous boundary condition enforcement functions \eqref{eq:phi_i} $\phi_i(x)$
    and \eqref{eq:psi_i} $\psi_i(x)$ for the interval $\Omega_i = (0,1)$.}
    \label{fig:psi}
\end{figure}

The additional scaling functions $\phi_i(x)$ and $\psi_i(x)$, also for the interval $\Omega_i = (0,1)$, are shown in Figure~\ref{fig:psi}.  We emphasize that, like the function $v_i(x)$, the inhomogeneous boundary condition functions $\phi_i(x)$ and $\psi_i(x)$ are not unique.  While exploring alternate choices for these functions goes beyond the scope of the present paper, we remark that our preliminary numerical experiments suggested that a smoothly varying $\phi_i(x)$ and $\psi_i(x)$ function is in general more effective than a Dirac delta function, as it minimizes the presence of sharp gradients, which can hinder the convergence of a NN.

\subsection{Incorporation of data into a PINN} \label{WDS:sec_data}

While the appeal of PINNs is at least partly that solution data are not required for training, we remark that it is possible to incorporate available data into the NN by adding a so-called data loss term $\mathcal{L}_d(\theta)$, as in~\eqref{eq:loss}. This data loss term takes the form of 
\begin{equation} \label{eq:data_loss}
    \mathcal{L}_d(\theta) := \frac{1}{d}\sum_{i=1}^d\left(\hat{u}(x_i; \theta) - u(x_i) \right)^2,
\end{equation}
where $d$ is the number of collocation points used to evaluate the data loss, and $u(x_i)$ is the snapshot data at $x=x_i$.  While, as written, \eqref{eq:data_loss} compares the NN solution with exactly one snapshot of the solution, additional snapshots can be obtained by generating solutions to the governing BVP \eqref{eqn:residual}--\eqref{eqn:BC} for different parameter numbers (e.g., different $Pe$), different boundary conditions and/or different source terms.  



\section{The Schwarz alternating method for PINN-PINN and PINN-FOM coupling} \label{WDS:sec:pinns_schwarz}

Our goal in this paper is to extend the Schwarz alternating method, described earlier in Section~\ref{WDS:sec:Schwarz}, to PINN-PINN and PINN-FOM couplings, following an overlapping domain decomposition of the underlying spatial domain $\Omega$.  Two coupling scenarios are possible: 
\begin{itemize}
    \item \textit{Coupling Scenario 1.} Use DD and the Schwarz alternating method to facilitate training of PINNs (offline).
    \item \textit{Coupling Scenario 2.} Use DD and the Schwarz alternating method to couple pre-trained subdomain-local PINNs (online).  
\end{itemize}
As discussed in Section~\ref{WDS:sec:intro}, Coupling Scenario 1 has been considered in several recent works, most notably~\cite{WDS:LiD3M, WDS:LiDeepDDM, WDS:Heinlein:2020}, whereas Coupling Scenario 2 was explored in~\cite{WDS:Wang:2022}. Herein, attention is focused on Coupling Scenario 1. Specifically, we explore the hypothesis that, by creating and coupling (using the Schwarz alternating method) several subdomain-local PINNs with each other and/or with FOMs, it may be possible to accelerate PINN training. For the targeted advection-diffusion problem~\eqref{eqn:residual}--\eqref{eqn:BC}, we are particularly interested in the advection-dominated (high P\'{e}clet) regime, in which it is well-known that PINNs are particularly difficult to train~\cite{WDS:Mojgani:2023}. 
We plan to explore Coupling Scenario 2 in a subsequent publication.  

In order to develop a PINN-PINN coupling algorithm using the Schwarz alternating method, we combine ideas presented earlier in Sections~\ref{WDS:sec:Schwarz} and~\ref{WDS:sec:pinns}. As in Section~\ref{WDS:sec_SDBC}, suppose the subdomain $\Omega = (0,1)$ has been partitioned into $n_D$ overlapping subdomains $\Omega_i$, so that $\Omega = \cup_{i=1}^{n_{D}} \Omega_i$, where $\Omega_{i} = (\gamma_{2i-2}, \gamma_{2i-1}) \in \mathbb{R}$ for $i = 1, ..., n_D$, with $\gamma_0 = 0$ and $\gamma_{n_{D}} = 1$. Let $\mathcal{L}_{r,i}$ and $\mathcal{L}_{d,i}$ denote the residual and data losses, defined earlier in~\eqref{eq:loss2} and~\eqref{eq:data_loss}, respectively, but restricted to subdomain $\Omega_i$. Let $NN_i(x; \theta)$ denote a NN within domain $\Omega_i$, and let $\hat{u}_i(x;\theta) \approx u_i(x)$ be the approximation of the solution $u(x)$ in $\Omega_i$. We provide some pseudo-code for our (offline) PINN training in Algorithm~\ref{alg:schwarz_pinn}.  The reader can observe that both the boundary loss term $\mathcal{L}_{b,i}(\theta)$ and the approximate PINN solution $\hat{u}_i(x)$ depend on the type of DBC enforcement selected, namely the variable $DBCtype$, which can take on three values: 
\begin{itemize}
    \item Weak DBCs (WDBC), corresponding to weakly-imposed DBCs on all boundaries (section \ref{WDS:sec_SDBC}),
    \item Strong DBCs (SDBC), corresponding to strongly-imposed DBCs on all system and subdomain boundaries (Section \ref{WDS:sec_SDBC}),
    \item Mixed DBCs (MDBC), corresponding to strongly-imposed DBCs at the system boundaries and weakly-imposed DBCs at the Schwarz boundaries (Section \ref{WDS:sec_SDBC} with $\phi_i(x) = \psi_i(x) = 0$).
    \end{itemize}

\begin{small}
\begin{algorithm}[h!]
    \begin{small}
    \SetAlgoLined
    {\bf Inputs:} overlapping DD of $\Omega$ into $n_D$ overlapping subdomains $\Omega_i = (\gamma_{2i-2}, \gamma_{2i-1})$; relaxation parameters $\alpha_r, \alpha_d \in [0,1]$; $DBCtype \in \{\text{WDBC}, \text{MDBC}, \text{SDBC}\}$; Schwarz convergence criteria. \\
    \hrulefill \\ 
    Initialize $\hat{u}_i(\gamma_{2i-1}, \theta) = 0$ for $i=1, ..., n_D$. \\
    Set $NN_{-i}(\cdot; \cdot) = NN_{n_D + 1}(\cdot; \cdot) = 0$. \\
    \While{Schwarz method unconverged} {
        \For{$i = 1, ..., n_D$}{ 
            Train PINN $NN_i(x; \theta)$ in $\Omega_i$ with loss 
            $$\mathcal{L}_i(\theta) := \alpha_r \mathcal{L}_{r,i}(\theta)+ (1-\alpha_r) \mathcal{L}_{b,i}(\theta) + \alpha_d \mathcal{L}_{d,i}(\theta),$$
            where  
            $$
            \mathcal{L}_{b,i}(\theta) = \left\{ 
            \begin{array}{cl}
                \mathcal{L}_{b,i}^{sys}(\theta) + \mathcal{L}_{b,i}^{sch}(\theta), & \text{if } DBCtype = \text{WDBC}, \\
                \mathcal{L}_{b,i}^{sch}(\theta), & \text{if } DBCtype = \text{MDBC}, \\
                0, & \text{if } DBCtype = \text{SDBC}, \\
            \end{array} \right.
            $$
            with 
            $$
            \mathcal{L}_{b,i}^{sys}(\theta)= \left\{ 
            \begin{array}{cl}
                NN_i(0; \theta)^2, & \text{if } i = 1, \\
                NN_i(1; \theta)^2, & \text{if } i = n_D, \\
                0, & \text{otherwise},
            \end{array}\right.
            $$
            and
            $$
            \mathcal{L}_{b,i}^{sch}(\theta)= \left\{ 
            \begin{array}{cl}
                \left( NN_i(\gamma_{2i-1}; \theta) - \hat{u}_{i+1}(\gamma_{2i-1}, \theta) \right)^2, & \text{if } i = 1, \\
                \left( NN_i(\gamma_{2i-2}; \theta) - \hat{u}_{i-1}(\gamma_{2i-2}, \theta) \right)^2, & \text{if } i = n_D, \\
                \left( NN_i(\gamma_{2i-2}; \theta) - \hat{u}_{i-1}(\gamma_{2i-2}, \theta) \right)^2 +\\ \left( NN_i(\gamma_{2i-1}; \theta) - \hat{u}_{i+1}(\gamma_{2i-1}, \theta) \right)^2, & \text{otherwise}.
            \end{array}\right.
            $$
            \\
            Set 
            $$
            \hat{u}_i(x;\theta) =  \left\{
            \begin{array}{cl}
                NN_i(x; \theta) , & \text{if }DBCtype = \text{WDBC}, \\
                v_i(x) NN_i(x; \theta), & \text{if } DBCtype = \text{MDBC}, \\ 
                v_i(x) NN_i(x; \theta) + \phi_i(x) \hat{u}_{i-1}(\gamma_{2i-2};\theta) + \\
                \psi_i(x) \hat{u}_{i+1}(\gamma_{2i-1}; \theta), & \text{if } DBCtype = \text{SDBC},
            \end{array}\right.
            $$
            where $v_i(x)$, $\phi_i(x)$ and $\psi_i(x)$ are given by \eqref{eqn:v_x}, \eqref{eq:phi_i} and \eqref{eq:psi_i}, respectively. \ 
            \\
            Interpolate $\hat{u}_i(x, \theta)$ onto $x = \gamma_{2i}$. \\
            Check convergence criteria for Schwarz iteration. Exit if converged.\\
        }
    }
    \caption{Alternating overlapping Schwarz PINN-PINN coupling algorithm}
    \label{alg:schwarz_pinn}
    \end{small}
\end{algorithm}
\end{small}

It is straightforward to modify Algorithm~\ref{alg:schwarz_pinn} to the case where a FOM is employed 
in one or more subdomains $\Omega_i$. While the details are omitted here for the sake of brevity, 
some numerical results for a PINN-FOM coupling are presented and discussed in Section \ref{WDS:sec_PINN-FOM}.

\section{Numerical Results} \label{WDS:sec_numerical}

In this section, we perform some numerical experiments aimed at understanding to what extent DD combined with the Schwarz alternating method can assist with the training of PINNs for the target advection-diffusion BVP~\eqref{eqn:residual}--\eqref{eqn:BC}.  

In order to perform the studies, we created a Python code\footnote{This code is available on github at the following URL: \url{https://github.com/ikalash/Schwarz-4-Multiscale}.} which invoked the {\tt TensorFlow} library \cite{WDS:Tensorflow} for PINN training. All of the NNs evaluated use the same hyperparameters. Each network had a 1D input layer to receive the spatial data, $x$, two hidden layers each with 20 nodes per layer, and a 1D output layer for the approximations of $u(x)$, similar to the NN architecture shown in Figure~\ref{fig:PINN}. Each hidden layer employed the swish activation function, given by  
\begin{equation}
    \sigma(z) := \frac{z}{1+e^{-\mu z}},
\end{equation}
in~\eqref{eq:weights_biases1}--\eqref{eq:weights_biases3}, for $\mu = 1$. When calculating the loss, the values of the relaxation parameters in~\eqref{eq:loss} were $\alpha_r=0.25$, 
\begin{equation}
    \alpha_b = \left\{
    \begin{array}[h]{ll}
        1 - \alpha_r, & \text{if using WDBCs or MDBCs},\\
        0, & \text{if using SDBCs},
    \end{array}\right. \hspace{1cm} 
    \alpha_d = \left\{
    \begin{array}[h]{ll}
        1-\alpha_r, & \text{if using data loss},\\
        0, & \text{otherwise},
    \end{array}\right.
\end{equation}
in all cases. The input data were not normalized, as they were already on the interval $[0,1]$. We chose to use the Adam optimizer for all NNs, with a constant learning rate of $0.001$. The number of training epochs for each NN per Schwarz iteration was 1024. For all experiments, $M=1024$ quasi-random uniform collocation points were selected within the full domain $\Omega = (0,1)$, subdivided among the subdomains $\Omega_i$, as inputs to evaluate the residual loss~\eqref{eq:loss_simple}.  For experiments in which PINN-FOM couplings were evaluated (Section~\ref{WDS:sec_PINN-FOM}), the FOM was generated by discretizing the governing PDE~\eqref{eqn:residual} using a second-order accurate finite difference approach with a mesh resolution of $h = \frac{1}{1024}$. This finite difference solution was also considered the ground truth against which all PINN approximations were compared in relative error calculations. For experiments in which the data loss term~\eqref{eq:data_loss} was included in the NN minimization problem, all FOM representations of this system were created using a second-order accurate backward finite difference approach with 1024 evenly spaced points on the domain $x\in(0,1)$. These snapshot data were then transferred to each subdomain $\Omega_i$ by interpolating the FOM solution to the collocation points in $\Omega_i$. The convergence criteria required both the Schwarz relative error to drop below a tolerance of $\delta = 0.001$, and the $L_2$ relative error in the NN approximation with respect to the reference FOM solution to drop below a tolerance of $0.005$. We allowed each model a maximum of 100 Schwarz iterations before cutting off training and declaring the model non-converged for a given problem case. 

The study summarized herein is aimed at understanding the relative impact of the following PINN- and coupling-related parameters on both the accuracy and the efficiency of the PINN training process: 
\begin{itemize}
    \item the number of subdomains, $n_D$,
    \item the size of the overlap region,
    \item the type of DBC enforcement in the PINN (WDBC, MDBC or SDBC), and 
    \item the impact of including the data loss term $\mathcal{L}_d(\theta)$ \eqref{eq:data_loss} in the loss function being minimized.  
\end{itemize}

\subsection{PINN-PINN coupling: parameter sweep study} \label{WDS:sec_overlap}

We begin by investigating the impact of the size of the overlap region and the number of Schwarz subdomains $n_D$ on both the accuracy and the efficiency of the resulting PINN-PINN coupling.  To do this, we first perform a parameter sweep study in which we vary $n_D$ between 2 and 5, and the percentage overlap of the subdomains between 5-50\% in increments of 5\%. To determine the boundaries of each subdomain, we calculate the size of each subdomain, $S_D$, and the size of the overlapping regions, $S_O$, as functions of the number of subdomains and the percentage overlap,
\begin{equation} 
    S_D = \frac{1}{n_D(1 - p_{O}) + p_{O}}\, ,
\end{equation}
\begin{equation}
    S_O = p_{O}S_D\, ,
\end{equation}
where $p_O$ is the percentage overlap expressed as a value between 0 and 1. Knowing $S_D$ and $S_O$, the boundary points, $\gamma_{2i-2}$ and $\gamma_{2i-1}$, which define each subdomain, $\Omega_i$, can then be determined by
\begin{equation} \label{eq:gamma_i_formulat}
    \begin{array}[h]{l}
        \gamma_{2i-2} = (i-1)(S_D - S_O), \quad \text{for} \quad i=1,\dots,n_D,\\
        \gamma_{2i-1} = \gamma_{2i-2} + S_D\, .
    \end{array}
\end{equation}

We also investigate the impact of using SDBCs vs. MDBCs vs. WDBCs on the system and Schwarz boundaries, as well as the impact of including the data loss~\eqref{eq:data_loss} in the loss function being minimized. We consider two P\'{e}clet numbers, namely $Pe = 10$ and $Pe = 100$, in this study.
Higher P\'{e}clet number problems are considered later in Sections \ref{WDS:sec_higherPe} and \ref{WDS:sec_PINN-FOM}.  

\begin{figure}
    \centering
    \includegraphics[width=1\textwidth]{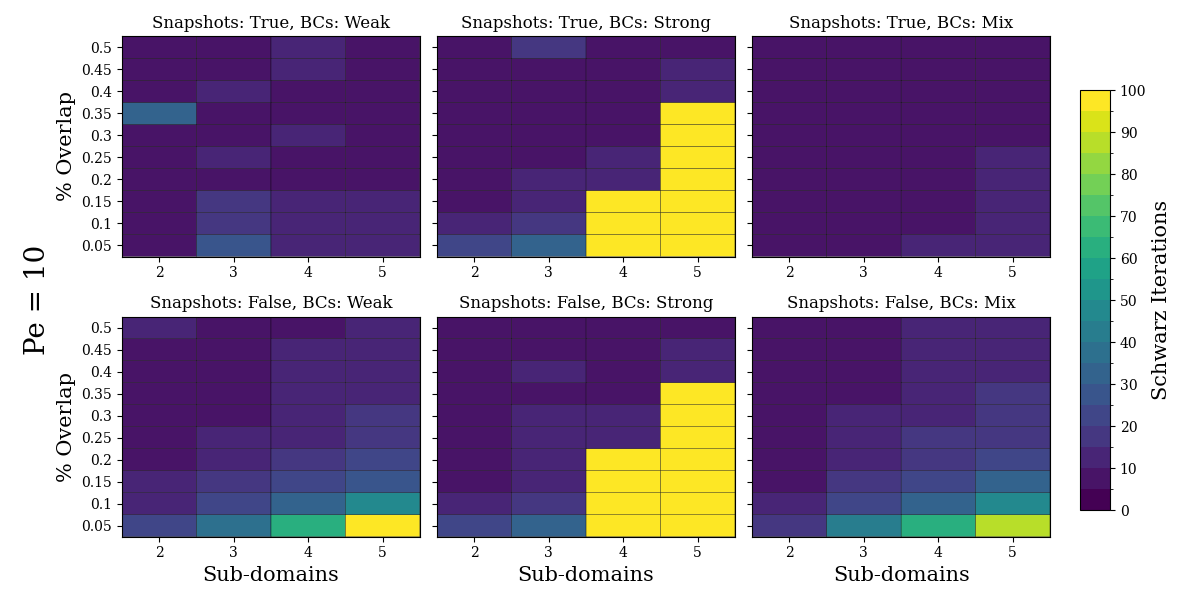}
    \caption{Schwarz iterations required to achieve convergence for various $n_D$ and percentage overlap, for $Pe = 10$. Each frame corresponds to various boundary conditions and data loss combinations.}
    \label{fig:iter_sweep10}
\end{figure}

\begin{figure}
    \centering
    \includegraphics[width=1\textwidth]{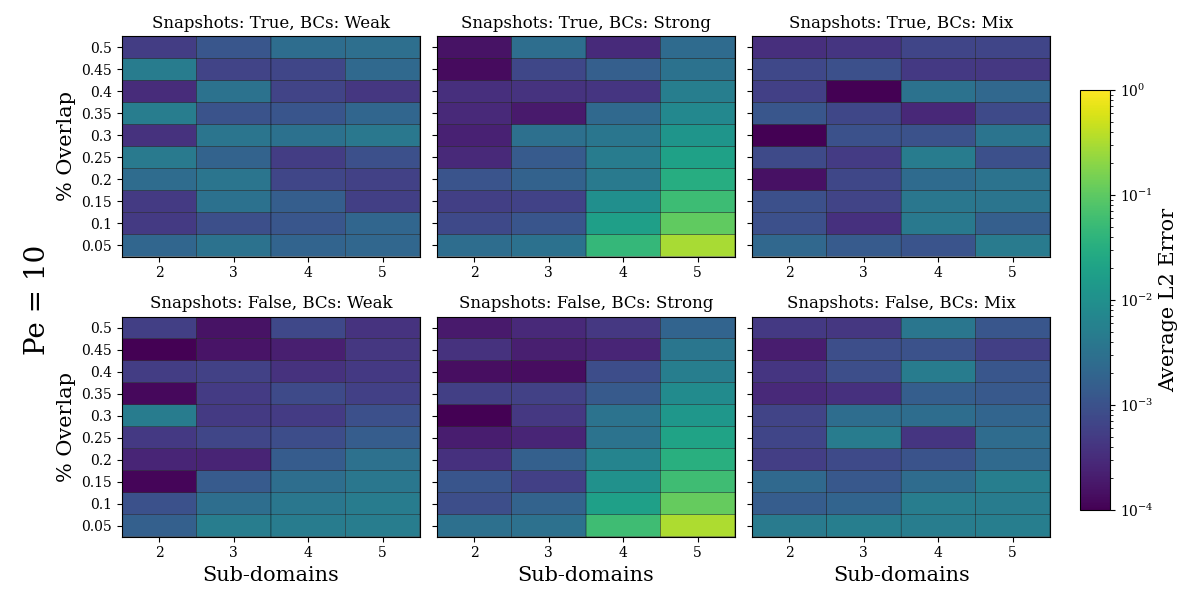}
    \caption{Average $L_2$ relative error at the end of the Schwarz iteration for various $n_D$ and percentage overlap, for $Pe = 10$. Each frame corresponds to various boundary conditions and data loss combinations.}
    \label{fig:l2_sweep10}
\end{figure}

\begin{figure}
    \centering
    \includegraphics[width=1\textwidth]{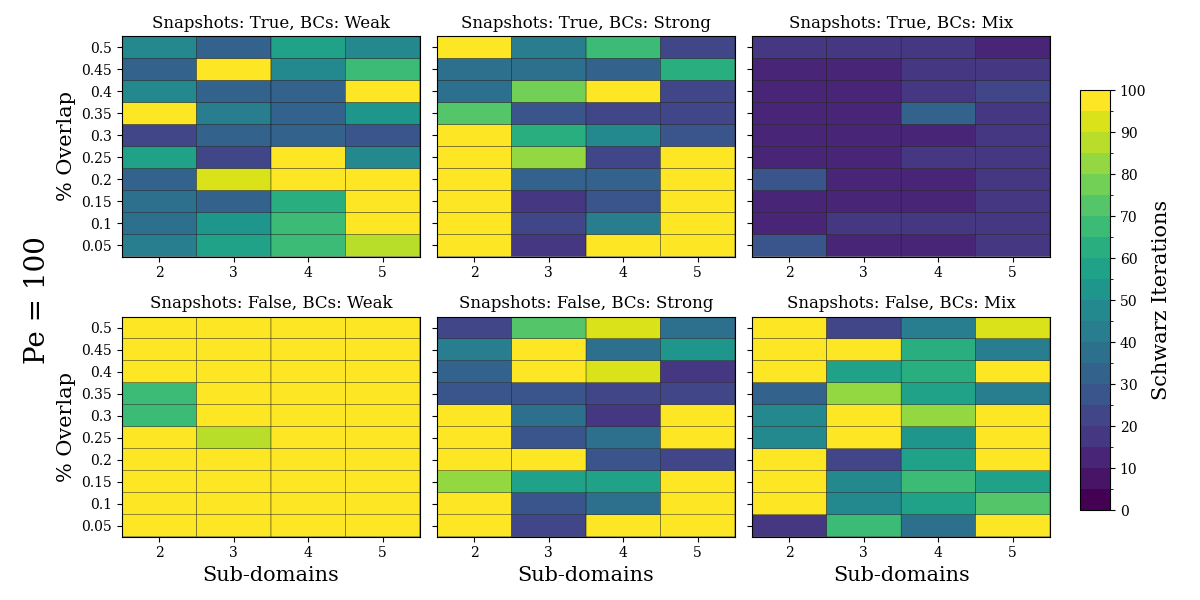}
    \caption{Schwarz iterations required to achieve convergence for various $n_D$ and percentage overlap, for $Pe = 100$. Each frame corresponds to various boundary conditions and data loss combinations.}
    \label{fig:iter_sweep100}
\end{figure}

\begin{figure}
    \centering
    \includegraphics[width=1\textwidth]{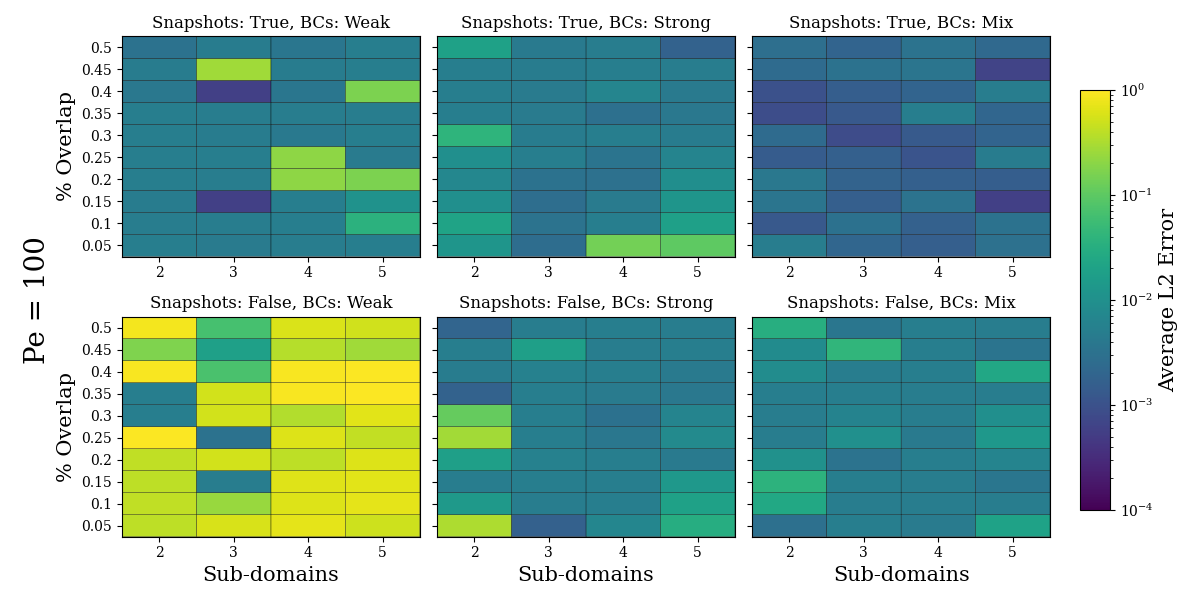}
    \caption{Average $L_2$ relative error at the end of the Schwarz iteration for various $n_D$ and percentage overlap, for $Pe = 100$. Each frame corresponds to various boundary conditions and data loss combinations.}
    \label{fig:l2_sweep100}
\end{figure}

To measure training efficiency and performance, both the number of Schwarz iterations required to satisfy the convergence criteria (Figures~\ref{fig:iter_sweep10} and~\ref{fig:iter_sweep100}) and the final $L_2$ relative error with respect to the FOM solution averaged across all subdomain models (Figures~\ref{fig:l2_sweep10} and~\ref{fig:l2_sweep100}) are recorded. For the case of $Pe = 10$, Figures~\ref{fig:iter_sweep10} and~\ref{fig:l2_sweep10} exhibit a general decrease in the number of iterations required for convergence and average $L_2$ relative error with increased percentage overlap. This is particularly true when SDBCs are enforced, as many cases with low overlap percentages fail to converge entirely (the cells marked in yellow in Figure \ref{fig:iter_sweep10}), but succeed with increased overlap. This follows the expected trend demonstrated in a number of references, including~\cite{WDS:Lions1988, WDS:mota2017schwarz, WDS:mota2022schwarz, WDS:LiDeepDDM}, wherein the number of Schwarz iterations (which typically correlates with the CPU time) is inversely proportional to the size of the overlap region. Conversely, increasing the number of subdomains seems to increase the required number of Schwarz iterations and the average $L_2$ error. This effect is, again, particularly clear for the case of SDBC enforcement. For the case of $Pe = 100$, however, Figures~\ref{fig:iter_sweep100} and and~\ref{fig:l2_sweep100} do not exhibit such clear relationships. Varying the number of subdomains and the percentage overlap appears to have a largely random effect, with many cases simply failing to converge within 100 Schwarz iterations (see the cells marked in yellow in Figure \ref{fig:iter_sweep100}). This is a stark example of the historic difficulty in constructing robust PINNs for fluid systems undergoing strong advection (see \cite{WDS:Mojgani:2023} and the references therein for more details).

Noticeable trends with respect to the boundary condition enforcement and data loss inclusion are also observed in our results. Curiously, SDBCs perform categorically worse than WDBCs or MDBCs in almost all cases. As noted previously, this is particularly true for the case of $Pe = 10$ with a larger number of subdomains and lower percentage overlap, for which many cases failed to converge. This conclusion generally extends to the case of $Pe = 100$, with the exception of WDBCs without a data loss term, for which nearly all cases failed to converge. Strikingly, the inclusion of a data loss term with WDBCs and MDBCs mostly eliminates this convergence problem. 
In general, for both $Pe = 10$ and $Pe = 100$, inclusion of the data loss term improves efficiency and accuracy. This is largely unsurprising, as one expects providing additional training data to improve the reconstructive accuracy of the neural network.
The fact that including the data loss term in the PINN loss function being minimized has the largest improvement for the models with WDBCs and MDBCs is consistent with the literature \cite{WDS:Faroughi2023, WDS:Rao2021}, which
reports that a PINN optimization problem is not guaranteed to be well-posed with a weak implementation of the boundary conditions.
We hypothesize that some of this ill-poseness is being encountered by our PINNs with WDBCs and MDBCs, and that the inclusion
of the data loss term helps to regularize the optimization problem in each of these cases.

\begin{figure}
    \centering
    \includegraphics[width=1\textwidth]{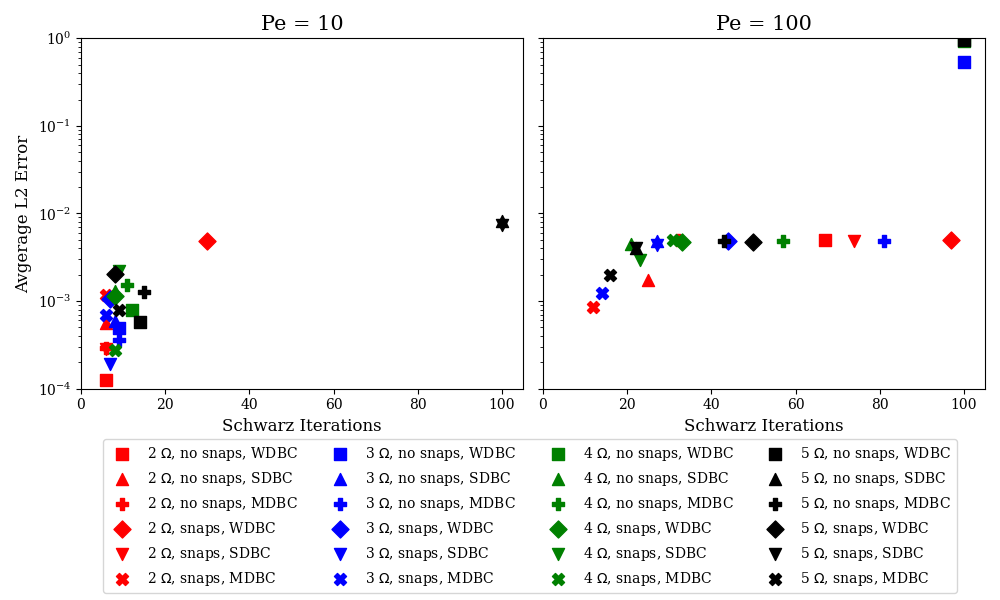}
    \caption{Pareto plots depicting $L_2$ relative error vs. Schwarz iterations for PINN-PINN couplings having different values of $n_D$ and different DBC enforcement types, with and without a data loss term, for $Pe = 10$ (left) and $Pe = 100$ (right).}
    \label{fig:pareto}
\end{figure}

To examine the above results in greater detail, Pareto plots for a single intermediate percentage overlap of 35\% display the average $L_2$ error with respect to the number of Schwarz iterations in Figure~\ref{fig:pareto}. The various symbol colors indicate the number of subdomains, while the symbol shapes indicate the boundary condition enforcement type and inclusion of the data loss. The apparent error threshold for $Pe = 100$ is an artifact of including the relative error in the convergence criteria. This broadly indicates the increased difficulty of efficient PINN training at higher $Pe$, where the predicted solution is slow to converge to the true solution, hence satisfying the Schwarz convergence criterion and yet struggling to satisfy the relative error criterion. These plots also draw attention to the relative success of the MDBCs in combination with the data loss term, indicated by an ``$\times$'' symbol. All such cases succeeded in converging in fewer than 40 iterations. 

\subsection{PINN-PINN coupling: the impact of different DBC implementations for higher $Pe$ numbers} \label{WDS:sec_higherPe}

Following the systematic parameter sweep studies described in the previous section, we performed a manual study exploring in more detail the impact of the boundary condition treatment on model convergence at high P\'{e}clet numbers. In this investigation, snapshot data are included in the loss function being minimized, and the relaxation term is set to $\alpha_r=0.2$. For these trials, the P\'{e}clet number is fixed at $Pe=150$, the number of subdomains is $n_D = 3$, and the percentage overlap is 20\%. Attempts to obtain a PINN-PINN coupled solution or a single-domain PINN solution for $Pe > 150$ were unsuccessful, as discussed in more detail in Section \ref{WDS:sec_singleDomainPINN}. As shown in Section \ref{WDS:sec_PINN-FOM}, it is possible to circumvent this difficulty by performing a PINN-FOM coupling using the Schwarz alternating method.  

\begin{figure}
    \begin{minipage}{0.49\linewidth}
        \includegraphics[width=0.99\linewidth]{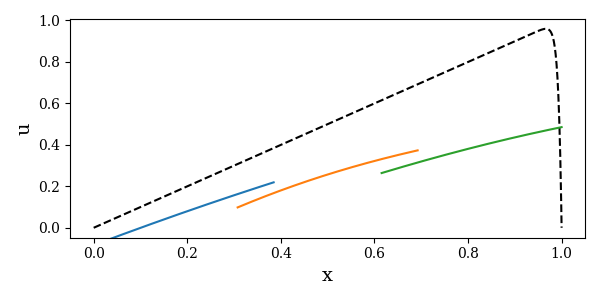}
        \subcaption{Iteration 1}
    \end{minipage}
    \begin{minipage}{0.49\linewidth}
        \includegraphics[width=0.99\linewidth]{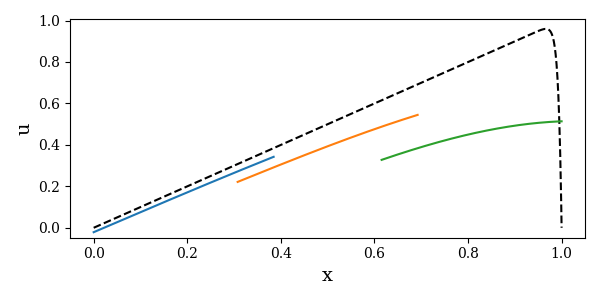}
        \subcaption{Iteration 25}
    \end{minipage}
    
    \begin{minipage}{0.49\linewidth}
        \includegraphics[width=0.99\linewidth]{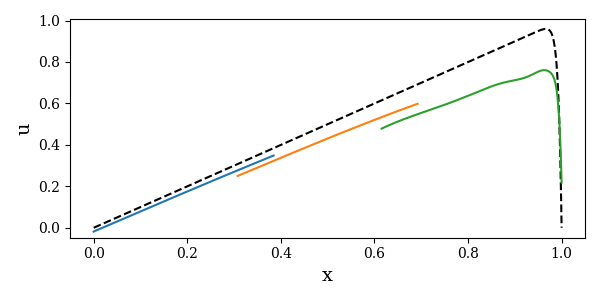}
        \subcaption{Iteration 50}
    \end{minipage}
    \begin{minipage}{0.49\linewidth}
        \includegraphics[width=0.99\linewidth]{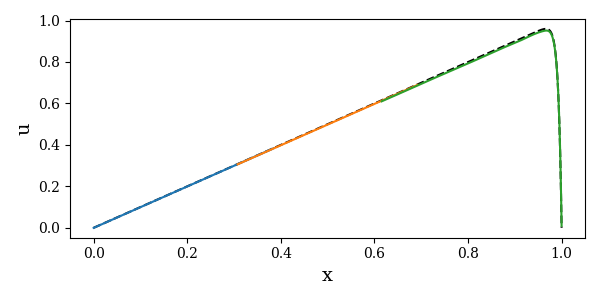}
        \subcaption{Iteration 70}
    \end{minipage}
    \caption{$Pe = 150$ alternating Schwarz-based converged solutions for a three subdomain all-PINN coupling with WDBCs for several Schwarz iterations.  The method converged in 70 Schwarz iterations.}
    \label{fig:WDBC_Pe150}
\end{figure}

\begin{figure}
    \begin{minipage}{0.49\linewidth}
        \includegraphics[width=0.99\linewidth]{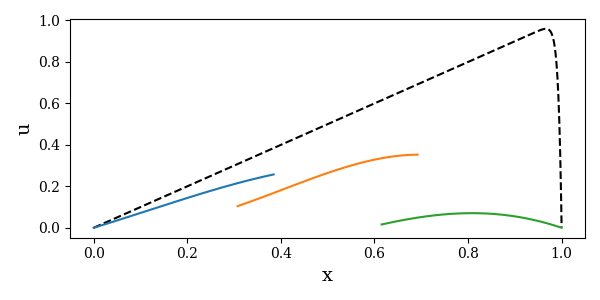}
        \subcaption{Iteration 1}
    \end{minipage}
    \begin{minipage}{0.49\linewidth}
        \includegraphics[width=0.99\linewidth]{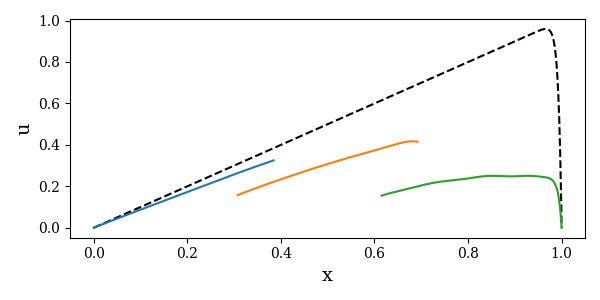}
        \subcaption{Iteration 10}
    \end{minipage}
    
    \begin{minipage}{0.49\linewidth}
        \includegraphics[width=0.99\linewidth]{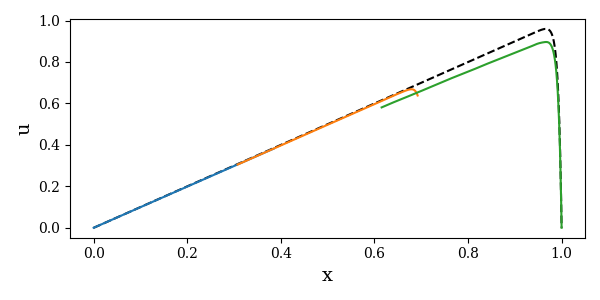}
        \subcaption{Iteration 20}
    \end{minipage}
    \begin{minipage}{0.49\linewidth}
        \includegraphics[width=0.99\linewidth]{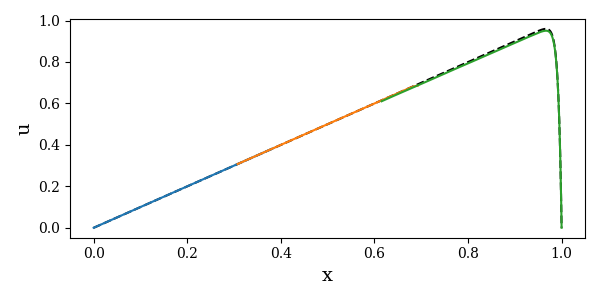}
        \subcaption{Iteration 26}
    \end{minipage}
    \caption{$Pe = 150$ alternating Schwarz-based converged solutions for a three subdomain all-PINN coupling with MDBCs (strong enforcement at system boundaries and weak enforcement at Schwarz boundaries) for several Schwarz iterations.  The method converged in 26 Schwarz iterations.}
    \label{fig:MDBC_Pe150}
\end{figure}

\begin{figure}
    \begin{minipage}{0.49\linewidth}
        \includegraphics[width=0.99\linewidth]{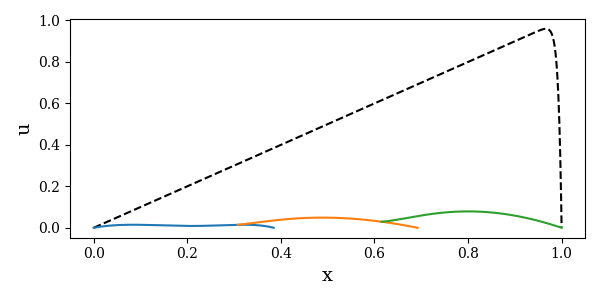}
        \subcaption{Iteration 1}
    \end{minipage}
    \begin{minipage}{0.49\linewidth}
        \includegraphics[width=0.99\linewidth]{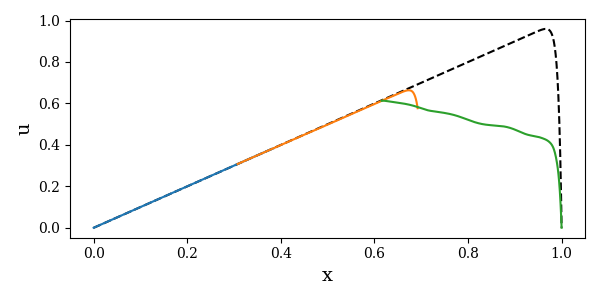}
        \subcaption{Iteration 25}
    \end{minipage}
    
    \begin{minipage}{0.49\linewidth}
        \includegraphics[width=0.99\linewidth]{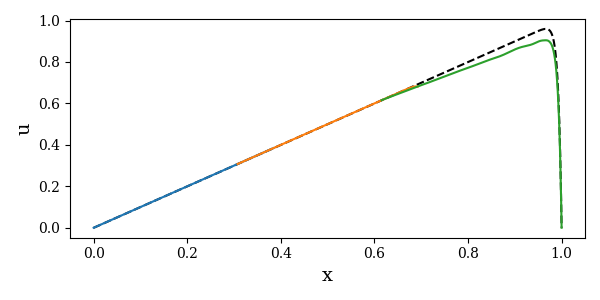}
        \subcaption{Iteration 50}
    \end{minipage}
    \begin{minipage}{0.49\linewidth}
        \includegraphics[width=0.99\linewidth]{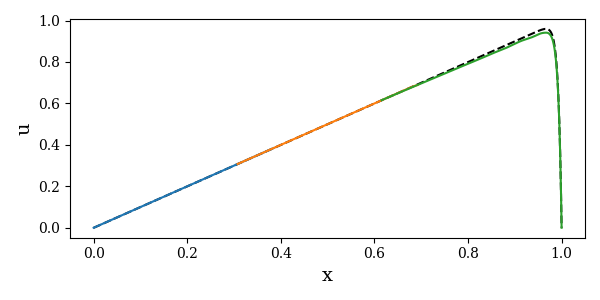}
        \subcaption{Iteration 68}
    \end{minipage}
    \caption{$Pe = 150$ alternating Schwarz-based converged solutions for a three subdomain all-PINN coupling with SDBCs for several Schwarz iterations.  The method converged in 68 Schwarz iterations.}
    \label{fig:SDBC_Pe150}
\end{figure}

Sample snapshots throughout model training are displayed for WDBCs in Figure~\ref{fig:WDBC_Pe150}, for MDBCs in Figure~\ref{fig:MDBC_Pe150}, and for SDBCs in Figure~\ref{fig:SDBC_Pe150}. The lower right panel in each figure represents the final solution upon achieving convergence. It is interesting to observe that, while all DBC enforcement strategies give rise to couplings that converged in fewer than 100 Schwarz iterations, the manner in which the solutions progress is very different. When imposing WDBCs, the subdomain-local PINNs very slowly converge to the true solution, and require additional iterations to closely match the system boundary conditions (Figure~\ref{fig:WDBC_Pe150}). In contrast, when imposing MDBCs, the subdomain-local PINNs very quickly approach a solution that resembles the FOM solution within the first 20 iterations, satisfying the system boundary conditions throughout convergence (Figure~\ref{fig:MDBC_Pe150}). Unlike the WDBC and MDBC cases, which display disjoint solutions between the subdomains, the enforcement of SDBCs (Figure~\ref{fig:SDBC_Pe150}) generates a relatively smooth solution throughout the duration of training. Similar to the MDBC case, the SDBC enforcement quickly approaches a reasonable solution, but requires a significant number of additional iterations to meet the convergence criteria.

\begin{figure}
    \begin{minipage}{0.49\linewidth}
        \centering
        \includegraphics[width=0.99\linewidth]{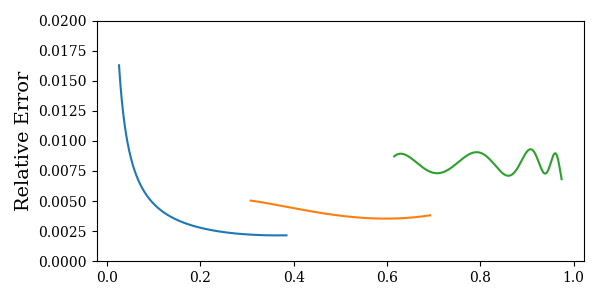}
        \subcaption{WDBCs, iteration 70}
    \end{minipage}
    \begin{minipage}{0.49\linewidth}
        \includegraphics[width=0.99\linewidth]{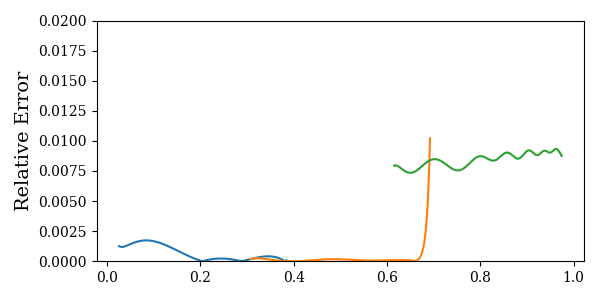}
        \subcaption{MDBCs, iteration 26}
    \end{minipage}

    \centering
    \begin{minipage}{0.49\linewidth}
        \includegraphics[width=0.99\linewidth]{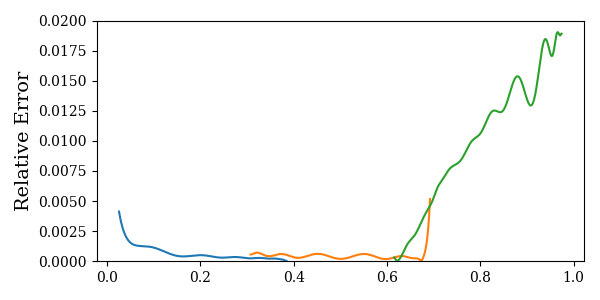}
        \subcaption{SDBCs, iteration 68}
    \end{minipage}
    \caption{Pointwise relative error in each subdomain for final $Pe = 150$ converged solutions displayed in Figs.~\ref{fig:WDBC_Pe150}, \ref{fig:MDBC_Pe150}, and \ref{fig:SDBC_Pe150}, in color and measured against the left $y$-axis. All couplings produced a solution with a relative error of $<2\%$ with respect to the exact solution.
    }
    \label{fig:pointwise_error}
\end{figure}

Pointwise relative error plots in Fig.~\ref{fig:pointwise_error} visualize the spatial accuracy induced by each boundary condition enforcement scheme. The error plots are trimmed to $x \in [0.025, \ 0.975]$ to avoid division by very small numbers in the error calculations. As one might expect, the WDBC and MDBC schemes do not guarantee that the error distribution is continuous between subdomains, as evidenced by significant jumps in relative error from one subdomain to the next. Conversely, the SDBC enforcement ensures that the predicted solutions match at the subdomain interfaces, resulting in a relatively smooth error profile. Also notable is the relatively high error near $x = 1$ when SDBCs are enforced, an artifact of the slow convergence of the solution in this region observed in Fig.~\ref{fig:SDBC_Pe150}. The pointwise relative error remains below 2\% in each case, and is ultimately controlled by the solution error tolerance used to determine convergence.

While all DBC enforcement methods investigated here lead to comparable results when considering Schwarz iteration counts alone, it can be argued that the convergence trajectory taken by the MDBC and SDBC variants are preferable. These models rapidly achieve a plausible solution, but struggle to fall strictly below the relative error convergence threshold. It is possible that a hybrid DBC enforcement strategy, in which one first applies MDBCs/SDBCs for a number of Schwarz iterations, followed by WDBCs in the subsequent Schwarz iterations, may give the most desirably convergence result, especially for higher P\'{e}clet numbers. 

\subsection{PINN-PINN coupling vs. single-domain PINN} \label{WDS:sec_singleDomainPINN}
One major question that has not yet been addressed is whether the proposed multi-domain coupling strategy for PINNs is actually beneficial for training such models in the advection-dominated regime when compared to a single-domain PINN having an equivalent architecture. To answer this question, we constructed a single-domain PINN for both forms of DBC enforcement and compared its convergence efficiency with the PINN-PINN couplings described above. Results were, however, largely inconclusive due to an observed extreme sensitivity to the model parameterization and initialization at high P\'{e}clet number, as previously displayed in Figure~\ref{fig:iter_sweep100}. We instead draw attention to broad trends. In general, strong enforcement of the system boundary conditions in the single-domain case resulted in faster convergence than those cases for which the boundary conditions were weakly enforced, in which case many models failed to converge entirely. For these high P\'{e}clet number systems, multi-domain PINNs with SDBCs similarly tended to converge faster than equivalent models with MDBCs or WDBCs, though they did not, on average, appreciably improve training efficiency over their single-domain SDBC counterparts. Although this appears to be a negative result for the proposed coupling strategy as a means to accelerate PINN training, we remark that it is likely possible to improve our PINN-PINN coupling results by introducing additional parallelism into the Schwarz iteration by employing the additive variant of the Schwarz alternating method (Remark 1), and/or alternating WDBC and SDBC enforcements during the training process, as suggested at the end of Section~\ref{WDS:sec_higherPe}. It may also be possible to attain improved results through PINN formulations especially designed for problems with a slowly decaying Kolmogorov $n$-width, such as the methods proposed in \cite{WDS:Mojgani:2023} and the references therein.

\subsection{PINN-FOM coupling} \label{WDS:sec_PINN-FOM}

\begin{figure}[!htbp]
    \centering
    \includegraphics[width=1\textwidth]{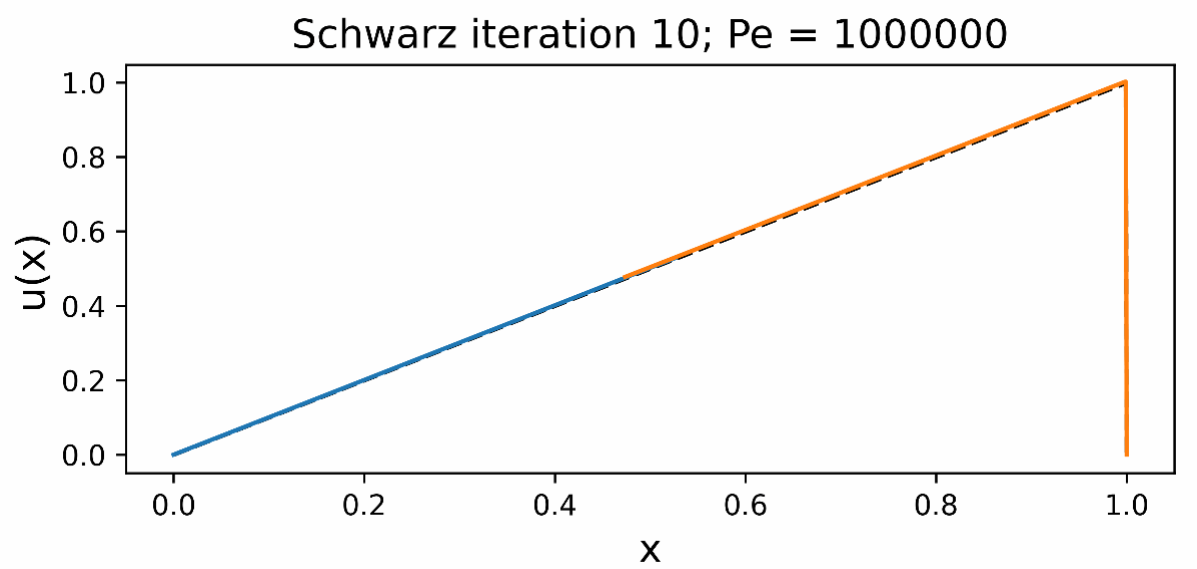} 
    \caption{Converged PINN-FOM coupled model with two subdomains and 10\% subdomain overlap for the advection-diffusion equation with $Pe = 1\times 10^6$. This model used WDBCs. The black dashed line indicates the FOM solution, the blue line indicates the PINN subdomain model solution, and the orange line indicates the FOM subdomain model solution.}
    \label{fig:pinn_fom}
\end{figure}

A coupling strategy enabled by the Schwarz alternating method but not yet investigated is the coupling of PINNs with high-fidelity FOMs.  
Our numerical results reveal for this coupling case that it is possible to obtain a convergent coupled model for arbitrarily-high P\'{e}clet numbers in as few as 10 Schwarz iterations by coupling a PINN to a sufficiently-accurate finite difference model used to represent the sharp gradient within the boundary layer that forms. Consider a DD of $\Omega = \Omega_1 \cup \Omega_2$ into $\Omega_1 = (0, \gamma_1)$ and $\Omega_2 = (\gamma_2, 1)$, where $\gamma_1$ and $\gamma_2$
are calculated using \eqref{eq:gamma_i_formulat} with $n_D = 2$ and $p_O = 0.1$. Figure~\ref{fig:pinn_fom} shows the result of a PINN-FOM coupling in which a PINN is specified in $\Omega_1$, and subsequently coupled to a finite difference model with mesh resolution $h = \frac{1}{1024}$ in $\Omega_2$. For this case, is possible to achieve convergence in just 10 Schwarz iterations for a P\'{e}clet number of $Pe = 1\times 10^6$.
We find that the type of DBC implementation (WDBC, MDBC or SDBC) has little effect on the accuracy and efficiency of the method.


The more interesting possibility that the PINN-FOM coupling strategy presents is the idea of pre-training several subdomain PINNs for different sections of a problem using FOM couplings on either side of the PINN. Thereafter, these pre-trained PINNs could be coupled online and may be capable of modeling more difficult problem cases than what is achievable with PINN-PINN coupled training, as discussed briefly in Section~\ref{WDS:sec_numerical}. We note, however, that this remains speculative, as we have not yet tested such a strategy. 

\section{Conclusions} \label{WDS:sec:conc}

In this paper, we explored the use of the Schwarz alternating method as a means to couple PINNs with each other and with conventional discretizations, following a decomposition of the spatial domain into overlapping subdomains.  The bulk of our attention was focused on determining whether domain decomposition and the Schwarz alternating method can expedite the offline training stage of PINNs, to facilitate training these models in the slowly-decaying Kolmogorov $n$-width regime (referred to as \textit{Coupling Scenario 1} herein).  
We explored three approaches for imposing Dirichlet boundary conditions within the overlapping Schwarz formulation at the system as well as the Schwarz boundaries: (i) a traditional weak formulation, in which a boundary loss term is added to the PINN loss function being minimized (termed WDBC), (ii) a strong formulation, motivated by the FBPINN literature \cite{WDS:Moseley:2023, WDS:Dolean:2023}, in which a transformation is applied to the neural network to ensure the solution satisfies the prescribed DBCs strogly (termed SDBC), and (iii) a ``mixed" formulation, in which the system DBCs are imposed strongly whereas the Schwarz DBCs are imposed weakly via the PINN loss function (termed MDBC).  We also explored the effect of incorporating a data loss term within the PINN loss function being minimized.  

We evaluated the proposed coupling approach on a canonical 1D advection-diffusion problem, in which sharp gradients that are a function of the P\'{e}clet number ($Pe$) form near the right boundary.  For $Pe = 10$ and $Pe = 100$, we performed a parameter sweep study in which we varied the number of subdomains, the size of the overlap region and the DBC implementation method.  The results indicate that at lower $Pe$, the domain decomposition has a strong effect on model training efficiency and accuracy, with larger overlap regions and a smaller number of domains tending to results in the fastest convergence to the true solution. This is largely consistent with the traditional Schwarz decomposition literature. At higher $Pe$, the decomposition structure has a less consistent effect. In almost all cases, however, including the data loss term resulted in a significant performance improvement. Curiously, the SDBC implementation of the boundary conditions did not always lead to the fastest convergence of the Schwarz alternating method, as hypothesized \textit{a priori}; the Pareto plot analysis reveals that the MDBC implementation leads to the best results in terms of accuracy and efficiency. While it is not clear from these results that Schwarz-based coupling of subdomain-local PINNs has an improvement on PINN training, especially at higher P\'{e}clet numbers, we demonstrate that PINN training \textit{can} be improved through a PINN-FOM coupling enabled by the Schwarz alternating method.  With this hybrid PINN-FOM coupling approach, it is possible to perform PINN training in as few as 10 Schwarz iterations for P\'{e}clet numbers as high as $1.0\times 10^6$.  

The initial exploration discussed in this paper has motivated several potential follow-on research efforts, itemized below.
 
\begin{itemize}
    \item \textit{Repeating the WDBC analysis with asymptotically-optimal norms and weights derived by Huynh, Meissner \textit{et al.} \cite{WDS:Meissner2023, WDS:Huynh2023}.} These weights and norms were derived by making analogies between PINNs and the Least Squares Finite Element Method \cite{WDS:LSFEMBook}.
    \item \textit{Extending the formulation to non-overlapping subdomains.}  The coupling of non-overlapping subdomains is important for multi-physics and multi-material problems, and requires the use of alternating Dirichlet-Neumann or Robin-Robin boundary conditions \cite{WDS:Gander2008, WDS:barnett2022schwarz}. 
    \item \textit{Extending the formulation to multi-dimensional and time-dependent problems.}  In multiple spatial dimensions, imposing Dirichlet BCs strongly is more challenging and will likely require a distance function-based approach \cite{WDS:Sukumar:2022, WDS:Wang:2022}.
    \item \textit{Prototyping of an additive Schwarz variant to improve efficiency.}  In the additive Schwarz implementation, it is possible to solve multiple subdomains in parallel, rather than contemporaneously, which can lead to greater CPU-time savings \cite{WDS:LiD3M, WDS:LiDeepDDM, WDS:Gander2008}.
    \item \textit{Considering alternate convergence criteria for the Schwarz alternating method.} Combining the Schwarz convergence and $L_2$ relative error in the convergence criteria may be sub-optimal, and alternative measures may be considered to improve predictive accuracy.
    \item \textit{Devising a ``hybrid" boundary condition enforcement strategy that switches between WDBCs and SDBCs.}  The results presented in this paper suggest that such a strategy has the potential to improve accuracy and efficiency.
    \item \textit{Exploring fully Coupling Scenario 2, in which pre-trained subdomain-local PINNs are coupled using the Schwarz alternating method.}  This work would leverage ideas from Wang \textit{et al.} \cite{WDS:Wang:2022}, in which a similar coupling scenario is investigated.
\end{itemize}

\section*{Acknowledgement} \label{WDS:sec:ackn}

Support for this work was received through Sandia National Laboratories' Laboratory Directed Research and Development (LDRD) program. The writing of this manuscript was funded in part by the second author’s (Irina
Tezaur’s) Presidential Early Career Award for Scientists and Engineers (PECASE). Sandia National Laboratories
is a multi-mission laboratory managed and operated by National Technology and Engineering Solutions of Sandia,
LLC., a wholly owned subsidiary of Honeywell International, Inc., for the U.S. Department of Energy’s National
Nuclear Security Administration under contract DE-NA0003525.
The authors wish to thank Alejandro Mota and Daria Koliesnikova for engaging in enlightening discussions related to details pertaining to the Schwarz alternating method for domain decomposition-based coupling, and Mamikon Gulian for offering invaluable insight into PINNs.  

\bibliographystyle{siam}
\bibliography{WillSnyder}


\end{document}